\documentclass[12pt]{article}
\usepackage{amssymb}
\usepackage{amsthm}\usepackage{epsf,epsfig,amsfonts}
\usepackage{amsmath}
\usepackage{graphicx}

\usepackage{a4wide}

\newtheorem{corollary*}{Corollary}

\newcommand{\be}{\begin{equation}}
\newcommand{\ee}{\end{equation}}

\newcommand{\weg}[1]{}

\usepackage{epsf,epsfig,amsfonts,a4wide}

\usepackage{amsfonts}
\usepackage{amsmath,amssymb,amsthm,color}
\usepackage{url}

\newtheorem{Th}{Theorem}

\newtheorem{Prob}{Problem}
\newtheorem{Subprob}{Subproblem}[Prob]

\theoremstyle{remark}
\newtheorem{Ex}{Example}
\newtheorem{Rem}{Remark}

\newcommand{\const}{\mbox{\rm const}}
\title{Geodesically equivalent metrics in  general relativity}
\date{} \author{  Vladimir  S. Matveev\thanks{ Institute of Mathematics, FSU Jena, 07737 Jena Germany,  vladimir.matveev@uni-jena.de}}
\begin{document}
\maketitle

\begin{abstract} 
We discuss  whether it is possible to reconstruct a metric by  its unparameterized geodesics, and how to do it effectively. We explain why this problem is interesting for general relativity. We show how to understand whether all curves from  a sufficiently big  family   are umparameterized geodesics of a certain affine connection, and how to reconstruct  algorithmically 
a      generic 4-dimensional metric by its unparameterized geodesics. The algorithm works most effectively 
if the metric is  Ricci-flat. We also prove that almost every  metric does not allow nontrivial geodesic equivalence, and  construct all   
pairs of 4-dimensional  geodesically equivalent metrics of Lorentz signature. 
\end{abstract} 

\section{Introduction}
\label{results} 

Let $(M^n,g)$ be a connected Riemannian  (= $g$ is positive definite)  or pseudo-Riemannian manifold of dimension $n\ge 2$.   We say that a metric $\bar g$ on $M^n$ is \emph{geodesically equivalent} to $g$, if every geodesic of $g$ is a (reparametrized) geodesic of $\bar g$.  We say that they are \emph{affine equivalent},  if their Levi-Civita connections coincide.

The first examples of 
geodesically equivalent metrics are due to  Lagrange \cite{lagrange}. 
He observed that the radial projection $f(x,y,z)= \left(-\frac{x}{z},- \frac{y}{z}, -1\right)$  takes geodesics of the half-sphere  
$S^2:=\{(x,y,z)\in \mathbb{R}^3: \ \ x^2+y^2+z^2=1, \ z<0\}$ to the geodesics of the plane $ E^2:=\{(x,y,z)\in \mathbb{R}^3: \ \ z=-1\}$,  since the geodesics of both metrics are intersection of the 2-plane containing the point $(0,0,0)$ with the surface.  Later, Beltrami \cite{Beltrami2,Beltrami3} generalized the example for the metrics of constant negative curvature, and for the pseudo-Riemannian metrics of constant curvature.    In the example of Lagrange, he  replaced the half sphere by  the half of 
 one of the hyperboloids $H_\pm^2:=\{(x,y,z)\in \mathbb{R}^3: \ \ x^2+y^2-z^2=\pm 1\}, $    with the restriction of the Lorentz metrics $dx^2+dy^2-dz^2$ to it. Then, the geodesics of the metric are also intersections of the 2-planes containing the point $(0,0,0)$ with the surface, and,  therefore, the stereographic projection sends it to the straight lines of the appropriate plane.
  
 Though   the examples of the Lagrange and Beltrami are   two-dimensional, one can easily generalize them for every dimension (for  Riemannian metrics, it was done already  in \cite{Beltrami2}) and for every signature.

Since the time of Hermann Weyl, geodesically equivalent metrics were actively discussed in the realm of geneal relativity theory. The context of general relativity poses the following restrictions: the dimension is $4$, the   metrics are pseudo-Riemannian of Lorentz signature $(-,+,+,+)$ or $(+,-,-,-)$, and  sometimes the metrics satisfy additional assumptions  such that one or both metrics are Ricci-flat ($R_{ij}=0$), or Einstein ($R_{ij}= \tfrac{R}{4}g_{ij}$), or, more generally, satisfy the Einstein equation $R_{ij}-\tfrac{R}{2}g_{ij} =T_{ij}$ with `physically interesting' stress-energy tensor $T_{ij}$.

Let us explain (using  a slightly  naive  language) 
 one of the possible motivations for this interest. Suppose we would like to understand the structure of the space-time  in a certain part of the universe. We assume that this part is far enough  so the we can use only telescopes (in particular we can not 
   send a space ship there). We still assume that the telescopes can see sufficiently many objects in this part of universe. Then, if  the relativistic effects are not negligible (that happens for example if 
   the objects in this 
  part of space time are sufficiently    fast 
   or  if this region of the universe is big enough), 
    we   obtain as a rule  the world lines  of the objects as unparameterized curves. 
    Indeed, local coordinates on a  4-manifold are  4  smooth functions on the manifold such that their differentials are linearly independent. 
    Now, for  every freely falling object in this part of the universe such that it can be registered by telescopes, each  telescope at every moment  of  time    
       gives  us two such functions,
        namely the spherical  coordinates $\phi$ and $\theta$ (latitude and longitude) 
    of the direction  the light reflected 
    from the object comes to the telescope  from  (in a naive language, the telescope `sees' the direction where the object lies), see the picture below. Since we have two telescopes, altogether we have   4 functions of $t$, \ $(\phi_1(t), \theta_1(t),\phi_2(t), \theta_2(t))$,  that we consider to be  
   the word line (i.e., geodesic) of the object in the coordinate system  $(\phi_1, \theta_1,\phi_2, \theta_2)$. If we see sufficiently many objects, we have sufficiently many geodesics. 
     
     Of course, we  cannot get lightlike or spacelike geodesics  by this procedure. In the best case, we can  reconstruct (numerically) sufficiently many geodesics, in the sense their velocity vectors are dense in a certain open subset of $TM$. See also the discussion in \cite{Hall2007}.

\begin{figure}
  \centering
  \includegraphics[width=\textwidth]{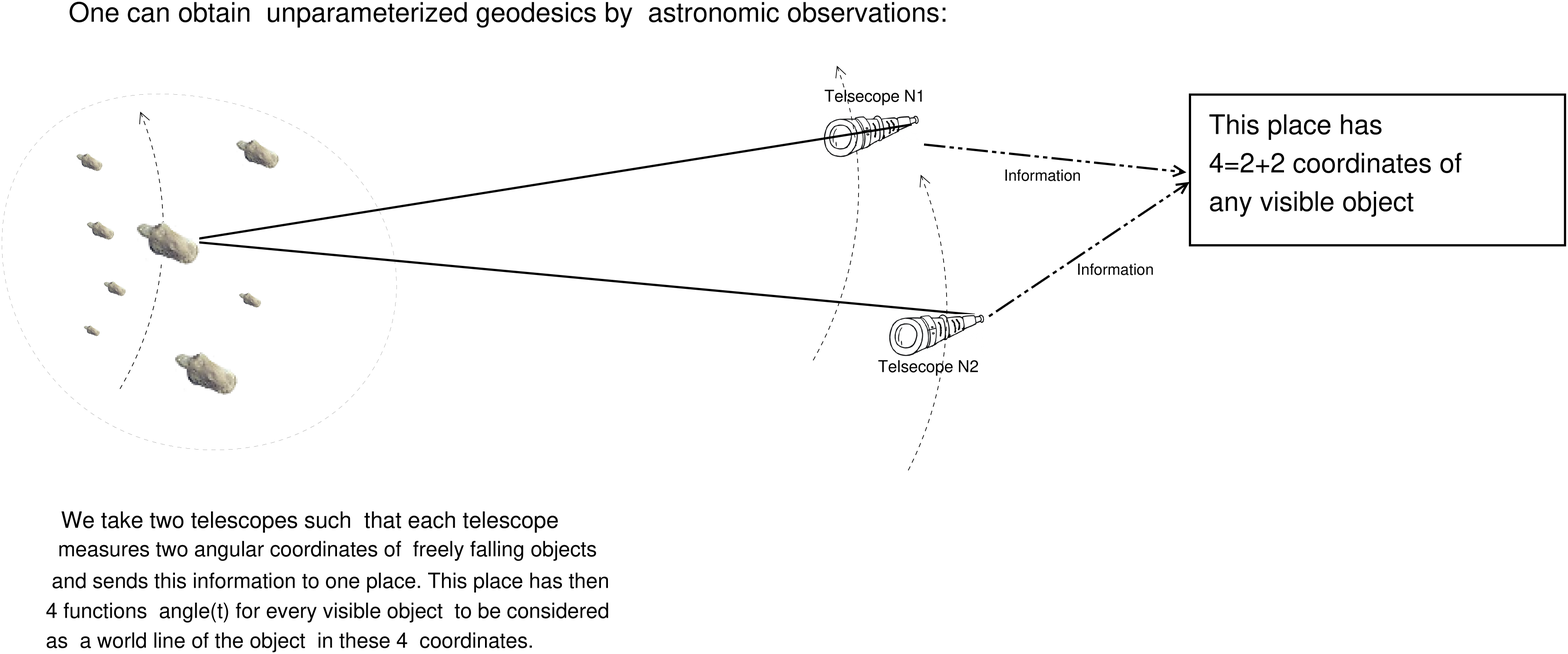}
\end{figure}
    
    Now, as a rule,  we can not get  the natural parameter (=proper time) of an object. Indeed,  if the relativistic effects are are not negligible, 
     the proper time of the object is not our own time $t$, i.e., the curve 
     $(\phi_1(t), \theta_1(t),\phi_2(t), \theta_2(t))$ is a reparameterized geodesic only. If we can not observe a periodic process on an object (note that  the astronomical  objects  such that we can  register a periodic process on, for example pulsars, are very rare)  or any other way to measure the own time of the object,
      we can not    obtain the own time of the objects by astronomic observations (see also the discussion in \cite{Gibbons}).

    In view of this discussion, the following two problems (Problem \ref{1} and Problem \ref{2} below)
     in the theory of geodesically equivalent metrics   are interesting for  general relativity: 
     
     \begin{Prob} \label{1} How to reconstruct a metric  by its unparameterized geodesics? 
     \end{Prob} 
     
     The general setting is as follows: we have   
   a family of smooth  curves $\gamma(t;\alpha)$ in $U\subseteq \mathbb{R}^4$ depending on 6-dimensional\footnote{locally, the set  of unparameterized 
    geodesics of an $n-$dimensional manifold has the structure of a manifold of dimension $2(n-1)$}
     parameter $\alpha=(\alpha_1,...,\alpha_6)$; 
     we assume that  the family is sufficiently big (we formalize `sufficiently big'  in the beginning of Section \ref{Problem11}). 
     We need to find a metric $g$ such that for every fixed 
      $\alpha$ the curve 
        $\gamma(t;\alpha)$ is a reparameterized geodesic of $g$.  
    
     Mathematically, the problem has sense in every dimension and for every signature of the metric.  
     In dimension 2, versions of this question were considered by S. Lie \cite{Lie} and R. Liouville \cite{Liouville}, and were  also discussed by Veblen and Thomas \cite{thomas,veblen23,veblen26} and   
      Eisenhart \cite{eisenhart23} in the beginning of 20th century. In the realm of general relativity, the problem was explicitly   stated by J. Ehlers et al \cite{Ehlers}, 
      where it was said  that   {\it ``We reject clocks as basic tools for setting up the space-time geometry   and propose ...  freely falling particles instead. We wish to show how the full space-time geometry can be synthesized ... . Not only the measurement of length but also that of time then appears as a derived operation.''}
      
      This problem can be naturally  divided  in two  subproblems. 
      
\begin{Subprob} \label{1.1} Given a family of curves $\gamma(t; \alpha)$, how to understand whether these curves are reparameterised geodesics of a certain affine connection? How to reconstruct this connection effectively? 
\end{Subprob}   

We will say that a metric  {\it lies in  a projective class} of a certain symmetric affine connection $\Gamma=\Gamma_{jk}^i$, if every geodesic of $g$ is a reparameterized geodesic of  $\Gamma$. 

\begin{Subprob}  \label{1.2}  Given  an affine connection $\Gamma=\Gamma_{jk}^i$, 
  how to understand whether there exists a metric $g$ in the projective class of $\Gamma$? How to reconstruct this metric effectively?
\end{Subprob}

Both subproblems were actively discussed in the literature.  In dimension 2,  the answer on  Subproblem \ref{1.1} is classical and was known already to Sophus Lie; given a family of  curves  one constructs an ODE  of the second order $y''(x) = f(x, y(x),  y'(x))$; the curves $\gamma(t; \alpha)$ are reparameterized geodesics of a certain connection if and only if the right hand side of the ODE is a 3rd degree polynomial  in $y'(x)$, $$f(x , y(x), y'(x))= A(x, y(x)) + B(x, y(x)) y'(x)+ C(x, y(x)) \left(y'(x)\right)^2+ D(x, y(x)) \left(y'(x)\right)^3.$$ 
The answer in the multidimensional  case  can be obtained  using  the same idea  as in dimension $2$, we give  it in Section \ref{Problem11}.

The second subproblem is  more complicated and is almost open.   In dimension 2, the  subproblem was considered  in the  
recent paper \cite{BDE} of Bryant et al: given an affine connection, 
       they construct a system of differential invariants 
      that vanish if and only if there exists a metric  (in a neighborhood of almost every point)  in the projective class of this connection.  The invariants are very complicated and are of  very high orders.

     In theory,  one can also obtain a similar answer in every dimension. Indeed,  by  \cite{eastwood}, 
      in every dimension the existence of a metric in a projective class is equivalent to the existence of a nontrivial solution of a certain overdetermined  system of linear PDE  in  the Cauchy-Frobenius form (i.e., the sysem is of first order and all  derivatives  of unknown functions 
      are explicit (linear) expressions in the unknown functions). Given an  overdetermined system of PDE in the  Cauchy-Frobenius form, one can always, in theory,  construct a system of differential invariants that vanish if and only if   the system admits a nontrivial solution (in a neighborhood of almost every point). An effective construction of these  differential invariants  could be very complicates. The results of \cite{BDE} show that it is indeed the case in dimension 2. It is hard  to predict whether the system of differential invariants is easier in the multidimensional case (normally multidimensional cases are  harder 
       than  lowdimensional;  
       but sometimes overdetermined systems are easier to analyse in 
       higher dimensions, because they can have higher degree of overdetermination). 
       
       In the present paper, in Section \ref{Problem12} we give an  algorithmic answer to Subproblem \ref{1.2}    under the additional assumption that the metric $g$ we are looking for is Ricci-flat and  the projective class  satisfies certain nondegeneracy assumption, i.e., in a situation most interesting from the viewpoint of general relativity. In Section \ref{Problem21}, 
       we also discuss the case of arbitrary metric: we show that  also in this case one can algorithmically reconstruct the metric by its projective class 
        assuming certain nondegeneracy assumption on  the projective class; though in this case the nondegeneracy assumption is harder to check.

       \begin{Rem} 
       Of course it is important in what form the geodesics $\gamma(t; \alpha)$ are  given. Below, it will be clear what information  we need   from $ \gamma(t; \alpha)$ in order our algorithm works. If the geodesics are given numerically (which is the case if they came from astronomic observations), this information could be extracted without difficulties.     
       \end{Rem}

       \begin{Prob} \label{2} In what situations is   the reconstruction of a  metric by the unparameterised geodesics unique (up to the multiplication of the metric by a constant)?  
     \end{Prob} 
        
        The example of Lagrange/Beltrami above shows that in certain situations the reconstruction is not unique:  the geodesics of  every 
        metric of constant curvature are  straight lines, i.e., the geodesics of the standard flat metric, 
        in  a certain coordinate system.
        Constant curvature metrics are not the only metrics that allow nontrivial 
         geodesical equivalence. For example, as it was shown by Dini, the following two metrics on $U^2\subseteq \mathbb{R}^2$  are geodesically equivalent    
     \begin{equation} \label{dini}  g= (X(x) - Y(y))(dx^2 +dy^2)   \ \textrm{and} \  
    \left(\frac{1}{Y(y)}  - \frac{1}{X(x)}\right)\left(\frac{dx^2}{X(x)} +\frac{dy^2}{Y(y)}\right),\end{equation} 
    where $X$ and $Y$ are arbitrary (smooth) functions of the indicated variables 
     such that the formulas \eqref{dini} correspond to metrics 
     (i.e., $0\ne X\ne Y \ne 0$  for all $(x,y)\in U^2$). 
  This  example was   generalized for all dimensions by Levi-Civita: from his results it is follows that the following two 4-dimensional metrics are geodesically equivalent:
     
 \begin{equation}\label{LC1}   \begin{array}{ccr} g &=&  (X_0(x_0)- X_1(x_1))(X_0(x_0)- X_2(x_2))(X_0(x_0)- X_3(x_3)) dx_0^2 \\ &+&  
(X_0(x_0)- X_1(x_1))(X_1(x_1)- X_2(x_2))( X_1(x_1)- X_3(x_3))dx_1^2 \\ &+& 
(X_0(x_0)- X_2(x_2))(X_1(x_1)- X_2(x_2))( X_2(x_2)- X_3(x_3))dx_2^2 \\ & + & (X_0(x_0)- X_3(x_3))(X_1(x_1)- X_3(x_3))( X_2(x_2)- X_3(x_3))dx_3^3 \end{array} \end{equation}
     
     \begin{equation} \label{LC2}      \begin{array}{ccr} \bar g &=&  \frac{1}{X_0(x_0)} \frac{1}{X_0(x_0)X_1(x_1)X_2(x_2)X_3(x_3)}  (X_0(x_0)- X_1(x_1))(X_0(x_0)- X_2(x_2))(X_0(x_0)- X_3(x_3)) dx_0^2 \\ &+&  \frac{1}{X_1(x_1)} \frac{1}{X_0(x_0)X_1(x_1)X_2(x_2)X_3(x_3)}
(X_0(x_0)- X_1(x_1))(X_1(x_1)- X_2(x_2))( X_1(x_1)- X_3(x_3))dx_1^2 \\ &+& \frac{1}{X_2(x_2)} \frac{1}{X_0(x_0)X_1(x_1)X_2(x_2)X_3(x_3)}
(X_0(x_0)- X_2(x_2))(X_1(x_1)- X_2(x_2))( X_2(x_2)- X_3(x_3))dx_2^2 \\ & + & \frac{1}{X_3(x_3)} \frac{1}{X_0(x_0)X_1(x_1)X_2(x_2)X_3(x_3)}(X_0(x_0)- X_3(x_3))(X_1(x_1)- X_3(x_3))( X_2(x_2)- X_3(x_3))dx_3^3. \end{array} \end{equation}
     Here $(x_0,...,x_3)$ are local coordinates and 
     the functions $X_i$ depend on the indicated variables and are such that the metrics have sense. 
     
     In view of this, 
     in the realm of general relativity,  Problem \ref{2}  can 
     be naturally  divided  in two  subproblems.

     We call  a metric $g$ {\it geodesically rigid,} 
     if  every metric $\bar g$,   geodesic equivalent to $g$, is proportional  to $g$.
     \begin{Subprob} \label{2.1} What  metrics  `interesting' for general relativity are { geodesically rigid}?   
     \end{Subprob} 
     \begin{Subprob}  \label{2.2} 
     Construct all pairs of nonproportional geodesically equivalent metrics. 
     \end{Subprob} 
     
     Let us comment on these subproblems. The part of the Supbproblem \ref{2.1} that is hard  or even impossible to formalize is the word ``interesting''.   Instead of formalizing this notion, let us  give few results in this direction. 
     
     Probably the metrics that are  most interesting in the context of general relativity are 
      Ricci-flat nonflat metrics. As it was shown by A. Z. Petrov in \cite{Petrov1} (see also \cite{Hall2007} and \cite{einstein}),   
      
\noindent{\it   4-dimensional Ricci-flat nonflat  metrics  of Lorentz signature can not be geodesically  equivalent, unless they are affinely  equivalent}

\noindent (two metrics are {\it affinely equivalent}, if their Levi-Civita connections coincide. Affine equivalent Ricci-flat 4-dimensional metrics are completely understood). 
       It is one of the results  Petrov  obtained in 1972 the Lenin   prize, 
       the most important  scientific award of the  Soviet Union,   for.

Recently, the answer of Petrov was generalized in \cite{einstein} (see also \cite{Hall2010}): it was shown that 

\noindent{\it if      $g$ and $\bar g$ are geodesically equivalent metrics on a $4-$dimensional manifold, and  $g$ is Einstein and of nonconstant curvature, then the metrics are affinely equivalent}.

\noindent Let us also give  an example of a  metric that is important for   general relativity and that  is not geodesically rigid.  This is the so-called Friedman-Lemaitre-Robertson-Walker metric 

\begin{equation} \label{RW1} g =  -dt^2 + R(t)^2 \frac{dx^2 + dy^2 + dz^2}{
1 + \tfrac{\kappa}{
4 } (x^2 + y^2 + z^2)} \ 
; \ \ \kappa = +1; 0;-1,\end{equation}
where  R = R(t) is a  real function (the scale factor)  of the `cosmic time' $t$. The metric is not geodesically rigid. Indeed, for every constant $c$ such that the formula below has sense, the metric \begin{equation} \label{RW2} \bar g =  \frac{-1}{(R(t)^2+c)^2} dt^2 + \frac{R(t)^2}{c(R(t)^2+c)} \frac{dx^2 + dy^2 + dz^2}{
1 + \tfrac{\kappa}{
4 } (x^2 + y^2 + z^2)} 
 \end{equation} is geodesically equivalent to $g$ (one can see it directly as it was done  for example \cite{Nurowski} or \cite{Hall2008}, see also discussion in \cite{Gibbons}.  Actually, the pair of geodesically equivalent metrics  (\ref{RW1},\ref{RW2}) is a  special case of geodesically equivalent metrics from Levi-Civita \cite{Levi-Civita}).

 For certain functions $R$, the metric \eqref{RW1} is  the main ingredient of the so-called  Standard Model of modern cosmology, and is of cause very interesting for general relativity.

The metrics  listed above, i.e., Einstein metrics and FLRW metrics,   are without any doubt interesting for general relativity. Of cause, there are other metrics that could be interesting for general relativity, and we consider  it very important 
 to understand what `interesting'  metrics are geodesically rigid. In the present paper, in Section \ref{Problem21}, 
 we prove that
   { \it  almost every 4-dimensional  metric is geodesically rigid.}
  
  \noindent Let us explain what we understand under almost every. Our  result is local, so we will work in a small neighborhood $U\subset \mathbb{R}^4$
  with fixed  coordinates $(x_1,...,x_4)$.  We consider a metric $g$ as the mapping $g:U\to \mathbb{R}^{\tfrac{n(n+1)}{2}}= \mathbb{R}^{10}$; the space $  \mathbb{R}^{\tfrac{n(n+1)}{2}}$ should be viewed as the space of symmetric $n\times n$-matrices. On the space of metrics (viewed as mappings) we consider the standard uniform $C^{2}-$topology: the metric  $g$ is  $\varepsilon-$close to the metric $\bar g$ in this topology, if  
  the components of $g$ and their first and second derivatives are  $\varepsilon-$close to that of $\bar g$.
   
 In the present paper,   we  prove that  
 
\noindent  {\it for any metric $g$  and  every   $\varepsilon >0$ there exists a metric $\hat g$ such that    $\hat g$   is $\varepsilon$-close to $g$ in the $C^2-$sense, and such that $\hat g$ is geodesically rigid.  Moreover,  there exists $\varepsilon'>0$ such that 
 every metric  that is $\varepsilon'-$ close to  $\hat g$ in the $C^2-$sense is also geodesically rigid. }

\noindent The result  is also true in dimensions $\ge 4$; the proof is essentially the same. 
 Now, concerning the lower dimensions, the result is true in dimension 3,
  if we replace the    uniform $C^2-$ topology by the  uniform $C^3$-topology.  The proof (will not be given here)   is based on the same idea. In dimension 2, the result is again 
  true,   if we replace the   uniform $C^2-$ topology by the   uniform $C^8$-topology. 
  
  This result was expected, at least if we replace $C^2-$topology by $C^\infty$-topology.  Indeed, by Sinjukov  \cite{sinjukov} and  Eastwood et al \cite{eastwood},     the existence of  a metric geodesically equivalent to a given one is equivalent to the existence of a nontrivial solution of a certain linear  system of partial 
  differential equations  in the  Cauchy-Frobenius  form  \eqref{prol},  whose coefficients are certain invariant 
  expressions in the components  of the given  metrics and their  derivatives.
   It is known that the existence of  the solution 
   of such system is equivalent to certain differential conditions on  coefficients, that is, on 
   the entries  of the metrics.  If there exists at least one  metric  that is geodesically rigid, 
    then the differential conditions are not identically fulfilled, and  
     almost every  (in the  $C^\infty-$ sense) metric  is  geodesically rigid.  
    Now, the existence of geodesically rigid metrics  in dimensions $n\ge 3$ is wellknown (at least since Sinjukov \cite{sinjukov54}).  The existence of geodesically equivalent metrics in dimension $n=2$ is more tricky; it follows from Kruglikov \cite{kruglikov2008} where all above mentioned differential conditions were constructed.  So in a  certain sense our result is the improving $C^\infty-$ closeness (which should be clear to experts, though we did not find a place where it is written)   to  $C^2-$closeness.

Let us now comment on Subproblem \ref{2.2}. First of all, the problem is very classical, and  was explicitly  asked by E. Beltrami\footnote{ Italian original from \cite{Beltrami}: 
La seconda $\dots$  generalizzazione $\dots$ del nostro problema,   vale a dire:   riportare i punti di una superficie sopra un'altra superficie in modo  che alle linee geodetiche della prima corrispondano linee geodetiche della seconda}  in  \cite{Beltrami}. 
In the Riemannian case, it was solved by Dini in dimension 2 and Levi-Civita in all dimensions. 
     More precisely, Dini has shown that  locally, in a neighborhood of almost every point of a two-dimensional manifold,
      every two  geodesically equivalent metrics are given by the form \eqref{dini} in a certain coordinate system. Levi-Civita has generalized this result to every dimension, we recall his result in Section \ref{Problem22}.

     Unfortunately, the  proofs  of Dini and Levi-Civita  require that the (1,1)-tensor $g^{i\ell} \bar g_{\ell j}$ is semi-simple (i.e., has no Jordan blocks), and that  all its eigenvalues are real.  If one of the metrics is Riemannian, this condition is fulfilled automatically. 
      Examples show  the existence of geodesically equivalent pseudo-Riemannian 
       metrics such that  the (1,1)-tensor $g^{i\ell} \bar g_{\ell j}$ is not semisimple or/and its eigenvalues are not real. The examples exist already in dimension 2: as it was shown\footnote{As is was explained  in \cite{pucacco}, essential part of the result could be attributed to  Darboux 
       \cite{Darboux}} in \cite{pucacco}, the  metrics from every  column  of the table

 \begin{tabular}{|c||c|c|c|}\hline &  \textrm{Liouville case} & \textrm{Complex-Liouville case} & \textrm{Jordan-block case}\\ \hline \hline
$g$ & $(X(x)-Y(y))(dx^2 -dy^2)$ &  $\Im(h)dxdy$ & $\left( 1+{x} Y'(y)\right)dxdy $
\\  \hline  $ \bar g$  &$ \left( \frac{1}{Y(y)}-\frac{1}{X(x)}\right) \left( \frac{dx^2}{X(x)} -  \frac{dy^2}{Y(y)} \right)$&
 \begin{minipage}{.3\textwidth}$-\left(\frac{\Im(h)}{\Im(h)^2 +\Re(h)^2}\right)^2dx^2 \\   +2\frac{\Re(h) \Im(h)}{   (\Im(h)^2 +\Re(h)^2)^2} dx dy  \\ +  \left(\frac{\Im(h)}{\Im(h)^2 +\Re(h)^2}\right)^2dy^2 $
 \end{minipage} &  \begin{minipage}{.3\textwidth}$  \frac{1+{x} Y'(y)}{Y(y)^4} \bigl(- 2Y(y) dxdy\\
    + (1+{x} Y'(y))dy^2\bigr)$\end{minipage}\\ \hline
\end{tabular}
      are geodesically equivalent  (we assume that the functions $X$ and $Y$ depend on the indicated variables only, and that the function $h$ is a holomorphic function of the complex variable $z=x +  i \cdot y$). Moreover, every pair  of 2-dimensional
       geodesically equivalent pseudo-Riemannian metrics  has this form in a neighborhood of almost every point in a certain coordinate system. 
    
    By direct calculations we see that the  (1,1)-tensor $g^{i\ell} \bar g_{\ell j}$ for these metrics 
 is semisimple with two real eigenvalues  in the Liouville case (we also see that the form of the metrics is very similar to \eqref{dini}, the only difference is the signature), has  two complex-conjugated eigenvalues in the Complex-Liouville case, and is not semisimple in the Jordan-block case. 
 
 Actually, certain authors consider that 
 the Subproblem \ref{2.2}  is also  solved; the solution is attributed to Aminova \cite{Aminova}. Unfortunaltely, the author of the present paper does  not understand her result, and has certain doubts that it is correct. 
  More precisely,   in view of \cite[Theorem 1.1]{Aminova}
  and  the formulas  \cite[(1.17),(1.18)]{Aminova} for $k=1$, $n=4$ and all $\varepsilon$s equal to $+1$, 
   the following two metrics $g$ and $\bar g$ given by the matrices  (where $\omega$ is an arbitrary function of the variable $x_4$).

$$ \left[ \begin {array}{cccc} 0&0&0&3\,x_{{3}}+3\,\omega \left( x_{{4}}
 \right) \\\noalign{\medskip}0&0&1&2\,x_{{2}}\\\noalign{\medskip}0&1&0
&x_{{1}}\\\noalign{\medskip}3\,x_{{3}}+3\,\omega \left( x_{{4}}
 \right) &2\,x_{{2}}&x_{{1}}&4\,x_{{1}}x_{{2}}\end {array} \right] ,
 $$

 {$
 \left[ \begin {array}{cccc} 0&0&0&3\,{\frac {x_{{3}}+\omega \left( x_
{{4}} \right) }{{x_{{4}}}^{5}}}\\\noalign{\medskip}0&0&2\,{x_{{4}}}^{-
5}&{\frac {-3\,x_{{3}}-3\,\omega \left( x_{{4}} \right) +2\,x_{{2}}x_{
{4}}}{{x_{{4}}}^{6}}}\\\noalign{\medskip}0&2\,{x_{{4}}}^{-5}&-{x_{{4}}
}^{-6}&{\frac {3\,x_{{3}}+3\,\omega \left( x_{{4}} \right) -2\,x_{{2}}
x_{{4}}+x_{{1}}{x_{{4}}}^{2}}{{x_{{4}}}^{7}}}\\\noalign{\medskip}3\,{
\frac {x_{{3}}+\omega \left( x_{{4}} \right) }{{x_{{4}}}^{5}}}&{\frac 
{-3\,x_{{3}}-3\,\omega \left( x_{{4}} \right) +2\,x_{{2}}x_{{4}}}{{x_{
{4}}}^{6}}}&{\frac {3\,x_{{3}}+3\,\omega \left( x_{{4}} \right) -2\,x_
{{2}}x_{{4}}+x_{{1}}{x_{{4}}}^{2}}{{x_{{4}}}^{7}}}&{\frac { \left( -3
\,x_{{3}}-3\,\omega \left( x_{{4}} \right) +2\,x_{{2}}x_{{4}} \right) 
 \left( 2\,x_{{1}}{x_{{4}}}^{2}+3\,x_{{3}}+3\,\omega \left( x_{{4}}
 \right) -2\,x_{{2}}x_{{4}} \right) }{{x_{{4}}}^{8}}}\end {array}
 \right]$}
  should be geodesically equivalent, though they are not (which can be checked by direct calculations).
  Note that the metrics above have signature $(2,2)$, so   they are not  that interesting for general relativity. In  the case of Lorentz signature, the theorem  of Aminova seems to be correct, but still it is very complicated to extract the precise formulas from her works. 

Note also that according to \cite{Aminova},  in the case of Lorentz signature, geodesically equivalent metrics were discribed by Petrov \cite{Petrov49} in dimension 3, by Golikov \cite{Golikov} in  dimension 4, and by Kruchkovich \cite{Kruchkovich} in all dimensions.  From these  papers, we were able to find (and to check) the paper of Petrov \cite{Petrov49}  only. 

{\it In the present paper,   we combine   recent 
 results   of \cite{splitting} and  above mentioned  results of \cite{pucacco} and  \cite{Petrov49} to 
  give an easy algorithm how to obtain  a list of pairs of all possible 
  geodesically equivalent 4-dimensional  metrics $g, \bar g$ of Lorentz signature. } 
  
  More precisely,  we explain (following \cite{splitting}) that every such pair can be obtained by applying the explicit 
   gluing construction from   Theorem \ref{thm3}  to building blocks, and provide  explicit formulas for all possible building blocks. One can easily obtain a complete list of metrics by this algorithm. There exists 
   three 
    possible three-dimensional building blocks, three possible two-dimensional, and one possible 1-dimensional, 
    so  all together there exists 10  normal  forms  for geodesically equivalent (nonproportional) 
    metrics of   Lorentz signature. The normal forms are given 
    by explicit formulas and allow certain freedom as (almost) arbitrary choice of functions of one variable or constants or metrics on  two- or three-dimensional disks.  We also explain  the (only) difficulty in applying this algorithm in higher dimensions.

    \section{ Problem \ref{1}: How to reconstruct  a metric by  its unparameterized geodesics. }\label{Problem1}

\subsection{ Subproblem \ref{1.1}: how to reconstruct a connection by unparameterized geodesics, and when it is possible. } \label{Problem11}

We  will work in arbitrary  dimension $n\ge2$, in a small neighborhood $U\subset \mathbb{R}^n$. 
We assume that we are  given a  family of smooth  
 curves $\gamma(t; \alpha)$.       We assume that  the family is sufficiently big in the sense that 
     at any point $x_0\in U$ the set of vectors 
     $$\Omega_{x_0}:= \{ \xi \in T_{x_0}U\mid \textrm{ there exists $\alpha$ and $t_0$ such that  $\tfrac{d}{dt} \big(\gamma(t;\alpha)\big)_{|t=t_0}  $ 
     is proportional to $\xi$}  \}$$ contains an open  
    subset of $T_{x_0}U$.  We put $\Omega= \bigcup_{x\in U}\Omega_x$.   We  will  call  a pair $(t_0; \alpha)$ \ {\it $x_0-$admissible}, if 
    $\tfrac{d}{dt} \big(\gamma(t;\alpha)\big)_{|t=t_0} \in \Omega_{x_0}$. 
          We need to understand whether there exists a symmetric 
    affine connection $\Gamma$ such that every curve $\gamma(t;\alpha)$ is a reparameterized geodesic of $\Gamma$, and construct this connection if it exists.

It is well known (at least since the time of Levi-Civita \cite{Levi-Civita}) that, 
in local coordinates, every  geodesic $\gamma:I\to U$,  $\gamma:t\mapsto \gamma^i(t) \in  U\subset  \mathbb{R}^n$   
 of a symmetric  affine connection  $\Gamma$  is given in terms of arbitrary parameter $t$ as  solution of 

\begin{equation} \label{arb} \frac{d^2  \gamma ^a}{dt^2}+ \Gamma_{bc}^a\frac{d\gamma^b}{dt}\frac{d\gamma^c}{dt}  = f\left(\frac{d\gamma}{dt}\right)\frac{d\gamma^a}{dt},\end{equation}  

  Better known version of this formula  assumes that the parameter is affine (we denote it by ``$s$'') and reads 
\begin{equation} \label{nat}  \frac{d^2\gamma^a}{ds^2}+ \Gamma_{bc}^a\frac{d\gamma^b}{ds}\frac{d\gamma^c}{ds}  = 0, \end{equation} 
it is easy to check that the change of the parameter $s \longrightarrow t$ transforms \eqref{nat} in   \eqref{arb}.

For further use, let us note that if we  linearly 
change the parameter $t$  of  a curve $\gamma(t; \alpha)$ (by putting $t= \const \cdot t_{new}$), 
the left hand side of  \eqref{arb} is   multiplied by $\const^2$ implying that the function $f$ should be  homogeneous of degree 1:
 $f(\const\cdot \xi)= \const\cdot f(\xi)$ for every $\xi$ (such that $\xi \in \Omega$). This allows us to 
  assume without loss of generality that for every $x$ the  
    subset  $\Omega_x\subseteq T_xU$   contains  a  cone over a nonempty  open subset.

  Let us now  take a point $x_0\in U$. For every $x_0-$admissible  $(t_0; \alpha)$, we    
  view the equations \eqref{arb}   as a system of equations  
    on the entries of $\Gamma(x_0)$ and on the function $f_{|\Omega_{x_0}}$; the coefficients in this system come from known data  $\left(\frac{d\gamma(t;\alpha)}{dt}\right)_{|t=t_0}$, $\left(\frac{d^2\gamma(t;\alpha)}{dt^2}\right)_{|t=t_0}$. Since we have infinitely many 
   $x_0-$ admissible  $(t;\alpha)$'s, we have an infinite system of  equations. Let us show that if  this system of equations is solvable, then  the solution is  unique  up to a certain `gauge' freedom. 
   
   Let us first describe the gauge freedom: we  consider two connections $\Gamma$ and $\bar \Gamma$ related by  Levi-Civita's formula 
  \begin{equation}  \label{bar} 
   \Gamma_{bc}^a= \bar \Gamma_{bc}^a - \delta_b^a\phi_c - \delta_c^a\phi_b,   
   \end{equation}  
    where   $\phi = \phi_i $ is 
a one form. Suppose the curve $\gamma$ satisfies  the equation  \eqref{arb} with a certain function $f$. 
Substituting  $\Gamma$ given by \eqref{bar} in    the left  hand side of   \eqref{arb}  and using
$$
(\delta_b^a\phi_c + \delta_c^a\phi_b)  \frac{d\gamma^b}{dt}\frac{d\gamma^c}{dt}= 
2 \left( \frac{d\gamma^b}{dt} \phi_b   \right) \frac{d\gamma^a}{dt} , 
$$ 
we obtain  that the same curve $\gamma$ satisfies the equation \eqref{arb} with respect to the connection $\bar \Gamma$ and  the function \begin{equation} \label{barf} 
\bar f(v):=  f(v) + 2 \left(v^b \phi_b   \right) .\end{equation} 
Thus, if $(\Gamma, f)$ is a solution of \eqref{arb}, then for every $1-$form $\phi$ 
the pair 
$\left(\bar \Gamma,  \bar f  \right)$  given by (\ref{bar},\ref{barf}) is also  a solution. 
 Let us show that up to this gauge freedom the connection $\Gamma$ and the function $f$  are  unique.

We again  work at one point $x_0\in U$  and 
 again view    \eqref{arb}  as equations on $(\Gamma,f)$. Suppose
 we have two  solutions $(\Gamma,f)$  and $(\bar \Gamma,\bar f)$. 
 We subtract  one equation from the other to obtain  
\begin{equation} \label{tilde}
\tilde \Gamma_{bc}^a v^b  v^c  = \tilde f(v)v^a, 
\end{equation}
where $\tilde \Gamma=  \bar \Gamma -\Gamma$, $\tilde f= \bar f- f$.  This equation is fulfilled for all 
 vectors $v= v^a$ lying  in an open nonempty $ \Omega_{x_0}\subset T_{x_0}U$.  Since the mapping 
$\sigma(u, v) \mapsto  \tilde \Gamma_{bc}^a u^b  v^c$  is linear in $u$ and $v$, it satisfies the  parallelogram
equality 
\begin{equation} \label{par} 0=\sigma(u +v, u+v) + \sigma(u -v, u-v) - 2\sigma(u, u) - 2\sigma(v,v).
\end{equation}
Combining \eqref{par}  with \eqref{tilde}, we obtain  
\begin{equation} \label{tildef1} \begin{array}{cl} 0&=\tilde f(u+v)(v+ u) + \tilde f(u -v)(u- v) - 2\tilde f(u) u - 2\tilde f(v) v \\ 
&=  (\tilde f(u+v)  +  \tilde f(u -v) - 2\tilde f(u))u + 
(\tilde f(u+v)  -  \tilde f(u -v) - 2\tilde f(v))v    .\end{array}\end{equation}  
Taking   $u$ and  $v$ to be  linearly independent, we obtain  
\begin{equation} \label{tildef2} \left\{\begin{array}{c} \tilde f(u+v)  +  \tilde f(u -v) - 2\tilde f(u)=0 \\  
\tilde f(u+v)  -  \tilde f(u -v) - 2\tilde f(v)=0\end{array} \right.  \end{equation} 
implying  $\tilde f(u+ v)= \tilde f(u)+ \tilde f(v)$. 
As we  explained above, 
 the functions $f, \bar f$, and, therefore, $\tilde f$,  also satisfy
  $\const \cdot \tilde f(v) = \tilde f(\const  \cdot v)$. Then,  
  the restriction of  $\tilde f$ to a certain nonempty open subset   $\Omega'_{x_0} \subset \Omega_{x_0}  \subset T_{x_0}U$ is linear, i.e., is given by 
  $\tilde f(v)= 2 \phi_a v^a$ for a certain 1-form $\phi=\phi_a$ and for all $v$ from $\Omega'_{x_0}$. 
   Then, the  connection 
  $$\hat \Gamma^a_{bc} := \bar \Gamma^{a}_{bc} - \phi_b \delta^a_c  - \phi_c \delta^a_b$$ 
  has the property  that for 
  every   $(t_0;\alpha)$  such that $\left(\frac{d  \gamma ^a}{dt}\right)_{|t=t_0}\in \Omega'_{x_0}$  the corresponding  $\gamma(t;\alpha)$ satisfies (at $t=t_0$) the equation 
  $$\frac{d^2  \gamma ^a}{dt^2}+ \Gamma_{bc}^a\frac{d\gamma^b}{dt}\frac{d\gamma^c}{dt}  = \frac{d^2  \gamma ^a}{dt^2}+ \hat \Gamma_{bc}^a\frac{d\gamma^b}{dt}\frac{d\gamma^c}{dt} \  \ \left(\Longleftrightarrow \   \Gamma_{bc}^a\frac{d\gamma^b}{dt}\frac{d\gamma^c}{dt}  = \hat \Gamma_{bc}^a\frac{d\gamma^b}{dt}\frac{d\gamma^c}{dt}\right) $$ implying $\Gamma = \hat \Gamma$ implying that $\Gamma$ and $\bar \Gamma$ are as in \eqref{bar} implying $f$ and $\bar f$ are as in \eqref{barf}.

  Finally, the connection $\Gamma$ and the function $f$, if they exist,     are uniquely determined by the unparameterized curves   $\gamma(t; \alpha)$ 
  up to the gauge freedom \begin{equation}\label{gauge}
  \Gamma_{bc}^a\mapsto   \Gamma_{bc}^a + \delta_b^a\phi_c + \delta_c^a\phi_b, \  \ f\mapsto f + 2 \phi\end{equation} 
  
  \begin{Rem} \label{f=phi} If the function $f$ is linear, i.e., if $f(\xi)= 
  2\phi_b \xi^b$  for a certain $1-$form $\phi$, then, up to the gauge freedom, we can take  $f\equiv 0$. 
  Moreover, putting  $f\equiv 0$ we exhaust the  gauge freedom.\end{Rem} 
  
  Let us now  explain how to  
  reconstruct the pair $(\Gamma, f)$ up to the gauge freedom.  We give an algorithm how to do it. 
  The algorithm  gives also a possibility to understand whether there exists such $(\Gamma, f)$: 
  we will see it that in order to uniquely reconstruct  the  (possible) 
  entries $\Gamma(x_0)^i_{jk}$ of the  connection  at a point $x_0$, we will need only finitely many $\gamma(t; \alpha)$ passing through this point.  There exists such $(\Gamma, f)$,  if for all $x_0$ the 
  entries of $\Gamma(x_0)^i_{jk}$ do not depend on the $x_0-$admissible    $(t_0; \alpha)$ we used to construct $\Gamma(x_0)^i_{jk}$.

  We   will work at 
   a point $x_0$; our goal is to reconstruct the components $\Gamma(x_0)_{jk}^i$. We take $x_0$-admissible 
   $(t_0; \alpha)$ such that  the first component 
   $\left(\tfrac{d\gamma^1}{dt}\right)_{|t=t_0} \ne 0$. For this geodesic $\gamma(t_0;\alpha)$,   
   we rewrite  the equation \eqref{arb}  at $t=t_0$ 
   in the  following form:  
  
  \begin{equation} \label{rec}
  \textrm{
  $ \begin{array}{rcl}  f\left(\tfrac{d\gamma}{dt}\right)  &=& \big(\tfrac{d^2\gamma^1}{d^2t} +   \Gamma_{ab}^1\tfrac{d\gamma^a}{dt}  \tfrac{d\gamma^b}{dt}  \big)/\tfrac{d\gamma^1}{dt} \\ 
 \tfrac{d\gamma^2}{dt}\,  \Gamma_{ab}^1\tfrac{d\gamma^a}{dt}  \tfrac{d\gamma^b}{dt} -  
 \tfrac{d\gamma^1}{dt}\,  \Gamma_{ab}^2\tfrac{d\gamma^a}{dt}  \tfrac{d\gamma^b}{dt} &=& \tfrac{d^2\gamma^2}{d^2t} \, \tfrac{d\gamma^1}{dt} - \tfrac{d\gamma^2}{dt}\,  \tfrac{d^2\gamma^1}{d^2t} \\ 
  &\vdots& \\ 
 \tfrac{d\gamma^n}{dt}\,  \Gamma_{ab}^1\tfrac{d\gamma^a}{dt}  \tfrac{d\gamma^b}{dt} -  
 \tfrac{d\gamma^1}{dt}\,  \Gamma_{ab}^n\tfrac{d\gamma^a}{dt}  \tfrac{d\gamma^b}{dt} &=& \tfrac{d^2\gamma^n}{d^2t} \, \tfrac{d\gamma^1}{dt} - \tfrac{d\gamma^n}{dt}\,  \tfrac{d^2\gamma^1}{d^2t}.  \end{array}$}\end{equation}
  The first equation of  \eqref{rec} is equivalent to the equation of \eqref{arb} for $a=1$
  solved  with respect to   $f\left(\tfrac{d\gamma}{dt}\right)$. We obtain the second, third, etc. equations of \eqref{rec} by substituting the first equation of \eqref{rec} in  the equations of \eqref{arb} corresponding to $a=2,3,\textrm{etc.}$
  
  We consider now a subsystem of \eqref{rec} containing the the second, third, etc. equations  of \eqref{rec}. 
  We see that the system  does  not contain the function $f$. Then, for every $x_0$-admissible $(t_0, \alpha)$,  it   is a  linear (inhomogeneous) system  on the components $\Gamma(x_0)_{jk}^i.$  
  We take a sufficiently big number $N$ and substitute $N$ \ 
    $x_0-$admissible  generic  $(t_0;\alpha)$'s  in this subsystem.

\begin{Rem}    If $n=4$, it is sufficient to take   $N= 12$. 
   We understand the world  {\it `generic}' in the  following sense: for every $n$   pairs $(t_0, \alpha)$, the velocity vectors    $\left(\tfrac{d\gamma}{dt}\right)_{|t=t_0}$ are linearly independent. \end{Rem}

    At every point $x_0$, we obtained  an  inhomogeneous linear system of equations on  $\frac{n^2(n+1)}{2}$ unknowns $\Gamma(x_0)^i_{jk}$. 
    
     {\it In the case  the  solution  of this system does not exist (at least at one point $x_0$),
      there exists  no connection 
 whose (reparameterized) geodesics are $\gamma(t; \alpha)$.     
     }
 
 If the solution exists at all points,  the solution is unique up to the gauge freedom \eqref{bar}. 
 Indeed, a solution of the last $n-1$ equations of \eqref{rec} gives us also the values  $f$ by  the first equation of \eqref{rec}, so the gauge freedom in the equations \eqref{rec} is the same as of the equations \eqref{arb}. Thus, a solution, if it exists,  
     gives us  the only up to the gauge  freedom   candidate for  the entries $\Gamma(x_0)_{jk}^i$  at every point $x_0$      such that its geodesics are (reparameterized)  curves $\gamma(t; \alpha)$.

    Assume now that at every point $x_0$,  a solution  $\Gamma(x_0)_{jk}^i$ exists. 
   In order to  construct the entries $\Gamma(x_0)_{jk}^i$ (up to the gauge freedome),
    we used $N$ \ $x_0-$admissible curves.  
  In order to understand whether all geodesics $\gamma(t; \alpha)$ are reparameterized geodesics of 
  $\Gamma$, we  need to substitute all 
   geodesics $\gamma(t; \alpha)$ in the equation \eqref{arb}, and check 
    whether it is  fulfilled; in this case, it is natural to rewrite the equation \eqref{arb} in the $f-$free  form $$\left(\frac{d^2  \gamma ^a}{dt^2}+ \Gamma_{bc}^a\frac{d\gamma^b}{dt}\frac{d\gamma^c}{dt}\right)   \wedge \frac{d\gamma^a}{dt}=0.$$

  \subsection{ Subproblem \ref{1.2}: 
 given  an affine connection $\Gamma=\Gamma_{jk}^i$, 
  how to understand whether there exists a metric $g$ in the projective class of $\Gamma$? How to reconstruct this metric effectively?} 
  
  \subsubsection{ General theory. }  \label{genth}
  We are given a  symmetric affine connection $\Gamma_{jk}^i$ on $M^n$, we need to understand whether there exists a metric in  the  projective class of $\Gamma$.     In this section we   recall (following \cite{bryant,eastwood}) the general approach how  to do it: the existence of a metric in  the projective class is equivalent to the existence of a nondegenerate 
   solution of a certain system of linear 
  PDE  in the Cauchy-Frobenius form, and, in theory,  there exists an algorithmic way to understand the existence of such solutions.

   \begin{Th}[\cite{eastwood}, see also references inside]  { $ g$ lies in a projective class of  
   a connection  {$\Gamma_{jk}^i$}  if and only if  
$\sigma^{ab}:= g^{ab} \cdot \det(g)^{1/(n+1)} $  is  a solution of }

\begin{equation} \label{east} {\left(\nabla_a\sigma^{bc}  \right) 
- \tfrac{1}{n+1}\left(\nabla_i\sigma^{i b}\delta^c_a + \nabla_i\sigma^{i c}\delta^b_a\right) =0.}  \end{equation}

Here $\sigma^{ab}:= { g^{ab}} \cdot { \det(g)^{1/(n+1)}}$ should be understood as an element of ${ S^2M } \otimes { (\Lambda_n)^{2/(n+1)}} M$.  In particular, 
$\nabla_a\sigma^{bc}= \underbrace{\frac{\partial }{\partial x^a} \sigma^{bc}  + \Gamma^{b}_{a{d}}\sigma^{{d}c} +  \Gamma^{c}_{{d}a}\sigma^{b{d}}}_{\textrm{\tiny   Usual covariant derivative}} - \underbrace{ \frac{2}{n+1}   \Gamma^{{d}}_{{d}a  }\, \sigma^{bc}}_{\textrm{ \tiny addition coming from volume  form}}$  
\end{Th}

    The equations \eqref{east} is a system of  $\left(\frac{n^2(n+1)}{2}- n\right)$ linear PDEs of the first order on $\frac{n(n+1)}{2}$ unknown components of  $\sigma$. 
    
    Two-dimensional  version of these equations was essentially known to R. Liouville \cite{Liouville}: instead of working with $\sigma^{ab}:= { g^{ab}} \cdot { \det(g)^{1/(n+1)}}$, he worked with $a_{ij}=\tfrac{1}{\det(\bar g)^{2/3}}\bar g_{ij}$; in dimension $2$ the entries of 
    $\sigma^{ij}$ and $a_{ij}$ are linearly related. The 2-dimensional analog of the equations \eqref{east} is then the Liouville  system of 4 PDE's of the first order 

\begin{equation}
\label{liouveq}   \begin{array}{rcc}
{\frac {\partial { a_{11}}}{\partial x{{}}}}  +2\,{{ K_0}}
{
a_{12}} -2/3\,{{ K_1}}  { a_{11}}  &=&0 \\
 2\,\frac
{\partial  a_{12}}{\partial x{{}}}  +\frac {\partial {
a_{11}} }{\partial y{{}}}  +2\,{ {K_0}}  { a_{22}} +2/3\,{
{K_1}} { a_{12}}  -4/3\,{{ K_2}}
  { a_{11}}  &=&0 \\
{\frac {\partial { a_{22}}}{\partial x{{}}}}  +2\,{\frac
{\partial { a_{12}}}{\partial y{{}}}} +4/3\,{ {K_1}}  { 
a_{22}}  -2/3\,{ {K_2}}  { a_{12}} -2\,{ {K_3}}
  { a_{11}}  &=&0 \\
{\frac {\partial { a_{22}}}{\partial y{{}}}}  +2/3\,{ {K_2}}
 { a_{22}}
 -2\,{{ K_3}}  { a_{12}}  &=&0,
\end{array}  
\end{equation} 
     where $K_0 :=-\Gamma^2_{11}$, $K_1:= \Gamma^1_{11}-2\Gamma^2_{12}$, $
K_2:= -\Gamma^2_{22}+2\Gamma^1_{12}$, $K_3:= \Gamma^1_{22}$.

\begin{Rem} 
 One sees that  the gauge freedom  \eqref{gauge} does not affect the coefficients $K_0,...,K_3$ of the equation  \eqref{liouveq}.  One can check by calculations that  this is  also true in all dimensions: the   
    gauge freedom \eqref{gauge}  does not change the equations \eqref{east}.  \end{Rem}

    The PDE-system \eqref{east}  can be prolonged (see \cite{eastwood}) to the system

\begin{equation} \label{prol} \left\{\begin{array}{rcl}
{\nabla_a} \sigma^{bc}&=&{ \delta_a{}^b}{\mu^c}
 +{ \delta_a{}^c}{\mu^b}\\
{\nabla_a}\mu^b&=&{ \delta_a{}^b}{\rho}
-{ \frac{1}{n}}{P_{ac}}\sigma^{bc}+{\frac{1}{n}}{W_{ac}{}^b{}_d}\sigma^{cd}\\
{\nabla_a}\rho&=&{ -\frac{2}{n}}{P_{ab}}\mu^b+
{\frac{4}{n}}{Y_{abc}}\sigma^{bc}
\end{array}\right.    \end{equation}

where  {$P$}  is the symmeterized   Ricci-tensor, {$Y$}  the  Cotton-York-Tensor and 
{$W_{ab}{}^c{}_d$} the projective   Weyl tensor for the connection   { $\Gamma$}. 

\begin{Rem} Here we use another index convention for the  projective   Weyl tensor than  in Section \ref{Problem12} of our paper. This convention is the same as in \cite{eastwood},   and is standard in the so-called tractor calculus, we refer to \cite{eastwood} for precise formulas. In Section \ref{Problem12} we will explain the   convention  used there 
by given the formula for Weyl tensor. 
\end{Rem}

The system \eqref{prol}  is a linear  system of PDE of the first order 
 on the unknown functions $\sigma^{bc}$, $\mu^b$ , $\rho$.  Moreover, 
 all  derivatives of unknowns  are expressed as functions of unknowns, i.e., the system is in the  Cauchy-Frobenius form.  
   One can understood this system geometrically 
    as a connection on the projective tractor bundle  ${\mathcal{E}}^{(BC)}
={\mathcal{E}}^{(bc)}(-2)+{\mathcal{E}}^b(-2)+{\mathcal{E}}(-2),$ see \cite{eastwood} for details. 
    The solutions of the system are then 
parallel sections of the connection; there exists an algorithmic way to understand whether a certain connection admits a  nontrivial parallel section. In the two-dimensional case, the algorithm was fulfilled for certain projectively homogeneous 
connections in \cite{bryant}; for arbitrary two-dimensional connection, the algorithm was fulfilled in \cite{BDE}, and the answer (i.e., the differential conditions on $K_i$ such that its vanishing implies the existence of a nontrivial solution) appears to be very complicated. 
In theory, one can fulfill this algorithm for every dimension; it is clearly a nontrivial task. In the next section we will show   that,   under the   additional assumption  that 
  the searched metric is Ricci-flat, there exists a trick that simplifies the algorithm.

 \subsubsection{ The case $n=4$, $g$ is Ricci-flat.} \label{Problem12} 
 Let us now assume that we know the geodesics of a  nonflat Ricci-flat metric. That is, we know a certain $\Gamma$ such that  for  a certain $\phi_a$ which we do not know 
 $\bar \Gamma_{bc}^a  :=  \Gamma_{bc}^a + \delta_b^a\phi_c + \delta_c^a\phi_c$ is the Levi-Civita connection of a certain nonflat 
 Ricci-flat metric which we again do not know.   Our goal is to find this metric (which I call $\bar g$).   By the above mentioned results of Petrov \cite{Petrov49}, Hall et al \cite{Hall2007,Hall2010},  and Kiosak et al \cite{einstein}, the metric  is unique up to multiplication by a constant; the goal of this section is to  explain how to find it algorithmically.  The algorithm  works under certain additional (generic)  condition on the connection $\Gamma$.
 
  We consider the projective Weyl tensor  introduced in \cite{Weyl} (not to be confused with the conformal Weyl tensor)

\begin{equation} \label{weyltensor} {W^i}_{jk\ell}:={R^i}_{jk\ell}-\tfrac{1}{n-1}\left({\delta^{i}_\ell} \, R_{jk}-{\delta^{i}_k}\, R_{j\ell}\right) \end{equation}  
(in our convention $R_{jk}=  {R^a}_{jka}$, so that ${W^a}_{jka}= 0$).

Weyl has shown that the projective Weyl tensor does not depend of the choice of connection within the projective class: if the connections $\Gamma $ and $\bar \Gamma$ are related by the formula \eqref{bar}, then their  projective Weyl tensors coincide.  
Now, from the formula \eqref{weyltensor}, we know that, if the searched $\bar g$ is Ricci-flat, projective Weyl tensor coincides with the Riemann tensor $ {\bar R^i}_{{ \  jk\ell}}$ of $\bar g$.    Thus, if we know the projective class of the Ricci-flat 
 metric $\bar g$, we know its Riemann tensor. 
 
 Then, the metric $\bar g$ must satisfy the following  system of 
 equations due to the symmetries of the Riemann tensor:

 \begin{equation} \label{equations} \left\{ \begin{array}{cc}  \bar g_{ia} {W^a}_{jkm} +  \bar g_{ja} {W^a}_{ikm} = 0   \\ 
 \bar g_{ia} {W^a}_{jkm} - \bar g_{ka}  {W^a}_{mij}=0\end{array}\right.  \end{equation} 
  
  The first portion of the eqautions \eqref{equations} is due to the symmetry 
 ($\bar R_{{ij}km}= -\bar R_{{ji}km}$), and the second portion is due to the symmetry  $(\bar R_{{km}ij}= \bar R_{ij{km}})$ of the curvature tensor of $\bar g$. 
 
 We see that for every point $x_0\in U$ \eqref{equations} is a system of linear equations on $\bar g(x_0)_{ij}$. 
The number of equations (around 100)  is much bigger than the number of unknowns (which is 10). It is expected therefore, that  a generic   projective Weyl tensor ${W^i}_{jkl}$ admits no more than one-dimensional space of solutions (by assumtions, our $W$ admits at least  one-dimensional space of solutions). The expectation is true, as the following classical  
result shows   
  
  \begin{Th}[\cite{Petrov1,Hall83,rendall,book,mcintosh}] 
  Let  ${W^i}_{jk\ell} $  be a tensor  in $\mathbb{R}^4$ such that it is skew-symmetric with respect to $k,l$ and 
  such that its traces ${W^a}_{ak\ell}$ and ${W^a}_{ja\ell}$ vanish. Assume that for all  1-forms   
  $\xi_i\ne 0$   we have ${W^a}_{jk\ell} \xi_a\ne 0.$ Then, the equations \eqref{equations} have no more than one-dimensional space of solutions. 
  \end{Th}

  Let us comment on the condition ${W^a}_{jk\ell} \xi_a\ne 0.$ In this context, for every fixed indexes $k,\ell$, ${W^i}_{j \ast \ast}$ could be viewed as a $n\times n$-matrix; and the condition ${W^a}_{j\ast \ast } \xi_a= 0$ means that the matrix has a nontrivial kernel (in particular, it is degenerate). Now, the condition ${W^a}_{jk\ell} \xi_a= 0$ means that for all indexes $k, \ell$ the kerns of the 
  $n\times n$-matrices ${W^a}_{jk\ell}$ have nontrivial intersection.  {\it Thus, it is a very restrictive condition on $W$, and, therefore, on $\Gamma$. }
   
 This  result shows that, under the assumptions that  for all 
  $\xi_i\ne 0$  we have  ${W^a}_{jk\ell} \xi_a\ne 0, $ we can reconstruct the conformal class of the metric $\bar g$ by solving the system of linear equations \eqref{equations}.  This can be done algorithmically. Then,  we also know the conformal 
  class of $\sigma$ 
  in  \eqref{east}, i.e., we know that $\sigma$ is of the form 
  \begin{equation} \label{ansatz} \sigma^{ij}= e^{\lambda} a^{ij},\end{equation} where $ a^{ij}$ is known and comes from the solution of the linear 
  system \eqref{equations}, and  the function 
  $\lambda$ is unknown. Substituting  the ansatz \eqref{ansatz} 
  in the  system   \eqref{east},  we obtain  an inhomogeneous system of linear equations on the components  $\tfrac{\partial\lambda}{\partial x^i}$. 
  Direct calculations show that this  system has at most one solution; since we assumed the existence of the metric in the projective class, one can always solve this  system  
and obtain   all  $\tfrac{\partial\lambda}{\partial x^i}$. Finally, we can obtain the function $\lambda$, and, therefore, the metric $\bar g$,  by integration. 
  
  Let us note that  in all steps we assumed that a  Ricci-flat 
  metric $\bar g$  exists   in the given projective class. 
  But the algorithm also gives us an algorithmic check whether such metric exists: one should go along the steps of the algorithm and look whether something goes wrong. 
  
  For example, the system \eqref{equations} could have no nontrivial 
  solution (i.e., every solution $\bar g_{ij} $ 
  of \eqref{equations} has zero determinant). Then,  no Ricci-flat 
  metric $\bar g$  exists in our 
   projective class.

  If the  system \eqref{equations} has nontrivial  
  solution, then, after plugging  the ansatz \eqref{ansatz} in \eqref{east}, we obtain a system of  nonhomogeneous linear 
  equations  on $\tfrac{\partial\lambda}{\partial x^i}$. This system may have no solution at all (the number of equations is much bigger than the number of unknowns; besides, the system is inhomogeneous), 
  or the $1-$form   $\tfrac{\partial\lambda}{\partial x^i}dx_i$
   may be  not  closed. In this case,  no Ricci-flat 
  metric $\bar g$  exists in our 
   projective class. 
   
   Finally, if the system \eqref{equations} has nontrivial  
  solution, if we can solve  the system of linear equaitons 
  we obtain after   plugging  the ansatz \eqref{ansatz} in \eqref{east}, and if the solution $\tfrac{\partial\lambda}{\partial x^i}$
   satisfies the `closeness' condition $ \tfrac{\partial }{\partial x^k} \tfrac{\partial\lambda}{\partial x^i}=\tfrac{\partial }{\partial x^i} \tfrac{\partial\lambda}{\partial x^k}$, then we do   obtain a 
    metric $g_{ij}$ in the projective class. The metric   must not be Ricci-flat though.  
   
   \begin{Rem} In Section \ref{Problem21},   we show that 
   one can reconstruct an almost every  (4-dimensional) 
   metric by its projective class, see Remark \ref{ssylka} there. In the case of arbitrary metric, the nondegeneracy assumption on the projective class is more complicated, and it is harder to check it.   
   \end{Rem}

 \section{ Problem 2: In what situations is the reconstruction of the metric by the unparameterised
geodesics unique (up to the multiplication of the metric by a constant)?}\label{Problem2} 

\subsection{ For generic 4-dimensional  metric,  the reconstruction of the metric by the unparameterized
geodesics is unique. }  \label{Problem21}
  
  Let us first construct one geodesically rigid metric in dimension $n=4$. 
  
  Using the formula \eqref{weyltensor}, by short tensor calculations we see that the metric $g_{ij}$  must satisfy the equation

\begin{equation}  \label{9} 
n g^{a (i}  W^{j)}_{\ \ \, a kl}   =g^{a b }W_{\ \ \, a b [l}^{(i} \delta_{k]}^{j)},
\end{equation} where $n=4$, the brackets ``$[ \ ]$'' denote the skew-symmetrization without division, and     the brackets ``$( \ )$'' denote the symmetrization without division.

\begin{Rem}
Actually, the equation (\ref{9}) is a part of the curvature of the tractor connection \eqref{prol}; in this context, it was obtained in \cite{eastwood}. 
\end{Rem}  

 We take a  4-dimensional  metric $\bar g$ such that     at the point $x_0$ it is given by the identity matrix 
 $$
    \begin{pmatrix} 1  &&&\\ &1&&\\ &&1&\\ &&&1\end{pmatrix}, 
   $$
   and such that its curvature tensor (with lowered indexes) $R_{ijkl} $
    at the point $x_0$ is given by 
 \begin{equation} \label{riem} 
R_{ijkl}=     h_{ik}h_{jl}- h_{il}h_{jk}+ H_{ik}H_{jl}- H_{il}H_{jk},
 \end{equation} 
    where  the entries  at $x_0$  of the  $(0,2)-$tensors $h$ and $H$ are  given by the  diagonal  matrices  
    $$
      h= \left[ \begin {array}{cccc} 1&&&\\ \noalign{\medskip}&2&&
\\ \noalign{\medskip}&&-1&\\ \noalign{\medskip}&&&0\end {array}
 \right] \ , \  
\ \ H= \left[ \begin {array}{cccc} 0&&&\\ \noalign{\medskip}&0&&
\\ \noalign{\medskip}&&1&\\ \noalign{\medskip}&&&1\end {array}
 \right]. 
    $$
    Such metric $\bar g$ exists  by \cite[Theorem 1.12.2]{Gilkey2001} 
    (see also \cite[Theorem 1.1]{Brozos}),  
    since  the  tensor \eqref{riem}  satisfies all symmetries of the curvature tensor.  
  
   Every metric $ g $ geodesically equivalent to $\bar g$ has the same projective 
  Weyl tensor as $\bar g$.  We view the equation \eqref{9} as the system of homogeneous linear equations  on the components of $g$; 
  every metric $g$ geodesically equivalent to  
     $\bar g$  satisfies this  system of equations (with the same coefficients  $W$!). 
     At the point $x_0$, 
     this  is a system on  $10$ unknowns $g(x_0)^{ij}$.  
   Since the system is symmetric in $i,j $ and 
  skew-symmetric in $k,l$, the system contains $60$ equations (actually, less because of certain hidden 
  symmetries inside).  By direct calculations, we see that the rank of this system is $9$. Indeed, it has at least one nontrivial  solution, namely $\bar g(x_0)^{ij}$,  so its rank is at most $9$. One can easily 
  find $9$ linear independent 
   equations of this system (so the rank is at least 9), namely the equations corresponding to the  followings  indexes $(i,j,k,l)$: 
   \begin{center} 
   \begin{tabular}{|c||c|}\hline  
                   $(i, j, k, l)$ & $\textrm{equation}$  \\\hline\hline 
                         $(1, 1, 2, 1)$ &  $-12 g^{1 2} = 0$  \\
                          $(1, 1, 3, 1)$ &  $2 g^{1 3} = 0$ \\ 
                          $(1, 1, 4, 1) $& $ 2 g^{1 4} = 0 $\\ 
                          $(2, 1, 2, 1) $ & $ 5 g^{1 1}-6 g^{2 2}+g^{3 3} = 0$\\
  $                         (2, 1, 4, 1)$ & $ g^{2  4} = 0 $\\ 
$                            (2, 2, 3, 2)$ & $ 8g^{2 3} = 0$ \\ 
                          $ (3, 1, 3, 1)$ & $-4g^{1 1}+g^{3 3}+4g^{2 2}-g^{4 4} = 0$\\
                          $ (3, 1, 4, 1) $& $ 2g^{3 4} = 0$ \\ 
                           $ (3, 2, 3, 2) $ &  $-6g^{2 2}+4g^{3 3}+3g^{1 1}-g^{4 4} = 0.$  
\\ \hline
\end{tabular} \end{center} 

   We see that the equations in the table  are linearly independent. 
   Thus, at the point $x_0$,   the set of solutions of this system is 1-dimensional, implying that 
    every metric $g$, geodesically equivalent to $\bar g$,  is proportional to $\bar g$. 
   
 Let us show that  at every point in a small neighborhood of $x_0$, the system \eqref{9} also has rank 9. Indeed, the rank of a matrix  is the biggest     dimension of  a nondegenerate quadratic submatrix and  therefore is  a lower semi-continuous (integer valued) function, i.e., rank  of this system is   
    at least 9 at every point of  a small neighborhood of $x_0$. Now, at every point the components  $\bar g^{ij}$ give us 
    a nontrivial  solution, so the rank can not be bigger than $9$.  Thus,  in a small neighborhood of $x_0$, 
    every metric $g$ geodesically equivalent to $\bar g$ is conformally equivalent  to $\bar  g$. Now, by 
     Weyl  \cite{Weyl}, two  conformally equivalent 4-dimensional metrics are proportional. Then, the metric $\bar g$ is geodesically rigid.

       Now let $\tilde g$ be an arbitrary metric in a small neighborhood of $x_0$. We consider the metric 
       $$
       g_{{t}}:= (1- t)\tilde g + t\bar g.  
       $$
        The system \eqref{9} constructed for this metric has rank $9$  for $t$ lying in a small interval  around 
         $1$.  Since the coefficients of the system are algebraic expressions in  $t$ whose 
         coefficients are algebraic expressions in 
         the components of $\bar g$, $\tilde g$ 
        and their first and second derivatives, for almost all $t$   the system \eqref{9} constructed for the 
        metric $g_{{t}}$ has rank $9$. We take   $t$  close to $0$ such that 
         the metric  $g_t$ is $\varepsilon-$close to $\tilde 
        g$ and  such that the system \eqref{9} constructed for the 
        metric $g_{{t}}$ has rank $9$. As we explained above, 
  this metric        is geodesically rigid.   Every metric $\hat g$ that is $C^2-$ close to $g_t$ is also geodesically rigid, since 
     the entries of $W$ for $\hat g$ 
       are algebraic expressions in the components of $\hat g$ and its first and second derivatives. Hence,       the 
      coefficients in the system   \eqref{9} constructed for $\hat g$  
       are close to that of   the system  \eqref{9} constructed for $g_t$ implying  the system 
        also has rank 9 implying    the metric $\hat g$ is geodesically rigid as well.
       
       {\it Thus, for every 4-dimensional metric  $\tilde g$   and for any $\varepsilon>0$ there exists a metric $g_t$ 
       that is $\varepsilon-$ close in the  $C^2-$sense to $\tilde g$ and    $\varepsilon'>0$ such that all metrics 
        $\varepsilon'-$ close in the  $C^2-$sense to $g_t$ are geodesically rigid. }

     \begin{Rem} 
     As we mentioned in the introduction, a similar proof can be done for all dimensions $n\ge 4$. For dimensions 2 and 3, the proof does not work anymore, since the system \eqref{9} has corank at least $2$ for all metrics $g$ (one can prove it using the methods of \cite[\S2.3.2]{KioMat2010}). One can still modify the proof replacing the system \eqref{9} by another projectively invariant 
     system of equations.  This other projectively invariant 
     system of equations requires higher derivatives of the components of $g$ though. 
      In dimension 3, one can construct (using the curvature of the tractor connection \eqref{prol}, see also \cite{Nurowski2})  such   
      projectively invariant system such that its coefficients 
       depend on the components of the metrics and its first, second and third derivatives. Therefore,  for 
       every  3-dimensional local  metric  $\tilde g$   and for any $\varepsilon>0$ there exists a metric $g_t$ 
       that is $\varepsilon-$ close in the  $C^3-$sense to $\tilde g$ and    $\varepsilon'>0$ such that all metrics 
        $\varepsilon'-$ close in the  $C^3-$sense to $g_t$ are geodesically rigid. Now, in dimension 2,     the construction of the projectively invariant system is much more involving (see \cite{BDE}) and requires 8 derivatives of the components of the metric. 
          \end{Rem}

  \begin{Rem} \label{ssylka}   We also see that the projective class of almost every (in the $C^2-$sense) 
   4-dimensional metric  determines its conformal class uniquely:  one can find the conformal class  by solving the system \eqref{9}.  Then, one can proceed along the 
     algorithm from Section \ref{Problem12}  and   understand whether there exists a metric in the projective class, and  find it.  
  \end{Rem} 
    
  \subsection{ Normal forms for  pairs of geodesically equivalent 4-dimensional  metrics such that one of them has Lorentz signature.}

   \subsubsection{ Splitting and gluing constructions  from \cite{splitting}.} \label{Problem22}
 Given two metrics $g$ and $\bar g$ on the same manifold, we consider the $(1,1)-$tensor $L=L(g,\bar g)$ defined by \begin{equation} \label{L}
L_j^i := \left(\frac{\det(\bar g)}{\det(g)}\right)^{\frac{1}{n+1}} \bar g^{ik}
 g_{kj},\end{equation}
 where ${\bar g}^{ik}$ is the contravariant inverse of
${\bar g}_{ik}$.

\begin{Rem}If $n$ is even, the tensor $L$ is always well defined. If $n$ is odd,  the ratio ${\det(\bar g)}/{\det(g)}$ may be negative, and  the formula \eqref{L} may have no sense. In this case, we replace   $\bar g$  by $-\bar g$ and  make the ratio ${\det(\bar g)}/{\det(g)}$ positive  and $L$  well defined.   In the cases interesting in  our context, $g$ and $\bar g$ have the same signature, and the problem with the sign does not appear at all.  \end{Rem} 

\begin{Rem} The tensor $L^i_j $  defined  in \eqref{L} is essentially the  same as 
 as the tensor  
introduced by Sinjukov (see equations (32, 34) on the  page 134 of the book \cite{sinjukov}, and also
Theorem 4 on page 135)  and which is often denoted by tensor $a_{ij}$ in the related literature. More precisely,  $L^i _j =a_{\ell j} g^{\ell i}$. It is also closely related to  $\sigma$ from  \S \ref{genth}:  $\bar g$ is geodesically equivalent to $g$, if and only if 
 $\bar \sigma^{ab} := {L^a}_{\ell}g^{\ell b}\cdot  \det(g)^{1/(n+1)}$ is a solution of  \eqref{east}.
 \end{Rem}

The  simplified  version of the  {\it gluing construction} does the following. 
Consider two manifolds $M_1$ and $M_2$ with  pairs of geodesically equivalent metrics $h_1\sim \bar h_1$ on $M_1$ and $h_2\sim \bar h_2$ on $M_2$. Assume  that the corresponding  $(1,1)$-tensor fields $L_1=L(h_1, \bar h_1)$ and $L_2=L(h_2, \bar h_2)$
 have no common eigenvalues in the sense that
for  any two points $x\in M_1$, $y\in M_2$ we have 
$$
\textrm{Spectrum}\, L_1(x) \cap \textrm{Spectrum}\, L_2(y) =\varnothing.
$$

Then one can naturally construct a pair of geodesically equivalent metrics $g\sim \bar g$ on the direct product   $M=M_1 \times  M_2$.  
These new metrics   $g$ and $\bar g$ differ from the direct product metrics  $h_1 + h_2$ and   $\bar h_1 + \bar h_2$ on $M_1\times M_2$ and are given  by the following   formulas involving $L_1$ and $L_2$: 
we denote by $\chi_i$, $i=1,2$, the characteristic polynomial of $L_i$: $\chi_i= \det(t\cdot {\bf 1} - L_i)$.   We treat the $(1,1)-$tensors $L_i$ as linear operators acting on $TM_i$.  
  A  polynomial $f(L)$ in $L$  is then the $(1,1)$-tensor of the form $f(L)=a_0(x) \cdot \mathrm{Id} + a_1(x) L + a_2(x) L^2 + \cdots + a_m(x) L^m$.
For two tangent vectors  $$u= (\underbrace{u_1}_{\in TM_1}, \underbrace{u_2}_{\in TM_2})\, , \  \ v=( \underbrace{v_1}_{\in TM_1}, \underbrace{v_2}_{\in TM_2}) \in TM $$
 we put    \begin{eqnarray} g(u,v) &  = &  h_1\left( \chi_2(L_1)( u_1), v_1\right)   + h_2\left(\chi_1(L_2)(u_2), v_2\right),  \label{hh1}  \\
  \bar g(u,v) &  = &  \frac{1}{\chi_2(0)}\bar h_1\left( \chi_2(L_1) (u_1), v_1\right)   + \frac{1}{\chi_1(0)}\bar h_2\left(\chi_1(L_2)(u_2), v_2\right).  \label{bh1}
\end{eqnarray}
  The corresponding  $(1,1)-$tensor 
$L=L(g,\bar g)$ is  the direct sum of $L_1$ and $L_2$  in the natural sense: 
for every 
$$
 v= (\underbrace{v_1}_{\in T_{x}M_1}, \underbrace{v_2}_{\in T_{y}M_2}){\in T_{(x,y)}(M_1\times  M_2)} \ \textrm{ \ we have\ } \ L(\xi)= \left(L_1(v_1), L_2(v_2)\right).$$

  It might   be  convenient to understand  the formulas (\ref{hh1}, \ref{bh1})  in matrix notation: we  consider   the coordinate system
$(x^1,...,x^r,y^{r+1},...,y^{n})$ on $M$  such that  $x-$coordinates are coordinates on    $M_1$ and  $y-$coordinates are coordinates on  $M_2$. 
Then, in this coordinate system, the matrices  of $g$ and $\bar g$  have the block diagonal form

 \begin{equation}\label{matg} 
g =\begin{pmatrix}   h_1   \chi_2(L_1) & 0 \\  0 &  h_2   \chi_1(L_2)\end{pmatrix}\ , \ \ \bar g =\begin{pmatrix} \frac{1}{\chi_2(0)}   \bar h_1 \chi_2(L_1) & 0 \\  0 & \frac{1}{\chi_1(0)}  \bar h_2   \chi_1(L_2)\end{pmatrix}.
\end{equation}

\begin{Th}[Gluing Lemma from \cite{splitting}] \label{thm3}

If $h_1$  is geodesically equivalent to $ \bar h_1$, and  $h_2$ is geodesically equivalent to $ \bar h_2$, then  the metrics $g,\bar g$ given by {\rm (\ref{hh1}, \ref{bh1})}  are geodesically equivalent too.
\end{Th}

 The {\it splitting construction} is the inverse operation.  We will not describe it completely  (and refer to  \cite{splitting}); we will use its  following corollary   explained in \cite[\S 2.1]{splitting}:
  
  {\it Every pair of 
 geodesically equivalent metrics    $h$ and $\bar h$ in a neighborhood of almost every point 
 can be obtained (up to a coordinate change) by applying  splitting construction  to  building blocks. } 
 
 By  a {\it building block} we understand an open neighborhood $U\subset \mathbb{R}^m  $ with a  pair of geodesically equivalent metrics $h\sim \bar h$   such that  at every point  the  tensor $L$ given by \eqref{L}  has only one real  eigenvalue, or two complex-conjugate
  eigenvalues, and such that the geometric multiplicity of the  eigenvalue  is constant  on  $U$. 

\begin{Rem} Riemannian version of the splitting/gluing constructions was known before, see  for example \cite[Lemma 2]{archive} and  \cite[\S\S2.2, 2.3]{bifurcations}. \end{Rem}

 \begin{Ex} \label{bb1} 
 In the definition of the building block, we allow the dimension $m=1$. 
 Then, the following two metrics on  the interval $I\subset \mathbb{R}^1$  with the following two 
 geodesically equivalent metrics $ h=dx^2$ and $\bar h  = X(x)dx^2$ (where  the function $X$ never vanishes) form a building block.    Actually, up to a coordinate change, 
 $(U_1, h, \bar h)$ is the only 1-dimensional building block. 
  \end{Ex} 
  
  \begin{Ex}  \label{bb2} 
All possible  examples of  two-dimensional  building blocks  can be extracted from the table of 2-dimensional geodesically equivalent metrics from the introduction. The metrics  from the first column  of the table do not correspond to a building block, since the tensor $L$ for these metrics has two different eigenvalues, $X(x)$ and $Y(y)$.  But the metrics from the second and the third columns do correspond to the building block, since the tensors $L$ for these metrics  are given by the matrices  
 $$
 \left[ \begin {array}{cc} {\Re(h)}&{ \Im(h)}\\ \noalign{\medskip}-{
 \Im(h)}&{ \Re(h)}\end {array} \right] \ , \ \ \left[ \begin {array}{cc} Y \left( x_{{2}} \right) &0
\\ \noalign{\medskip}1+x_{{1}}{\frac {d}{dx_{{2}}}}Y \left( x_{{2}}
 \right) &Y \left( x_{{2}} \right) \end {array} \right].
 $$
 
 Of cause,  in every dimension, in particular in dimension two,  there exists a trivial building block 
 $(U, h, \bar h = \const \cdot h)$; the  tensor $L$ for this metric is a multiple 
of   $\delta^i_j$.  
 From the results of \cite{pucacco} it follows that every  two-dimensional building block has one of these three forms.   
 \end{Ex}

 The formulas for the 3-dimensional building block can be obtained using  Petrov  \cite{Petrov49} and Eisenhart \cite{eisenhart38}; we will give them  later. From linear algebra it follows that if the metrics $g$, $\bar g$ have Lorentz signature, then  4-dimensional building blocks are not possible (except for the trivial block corresponding to proportional metrics $g\sim  \bar g:= \const \cdot g$), since in the Lorentz signature  a $g$-selfadjoint $(1,1)$tensor $L$ can not  have a Jordan block  of dimension $\ge 4$ with real eigenvalue, and a Jordan block of dimension $\ge 2$ with complex eigenvalue.

 \begin{Ex}[Dini formulas \eqref{dini} follow from splitting-gluing constructions.] \label{diniformulas} 
 We consider the two 1-dimensional building blocks
 $$
\left(I_1, h_1=dx^2, \bar h_1 = \frac{1}{X(x)^2}dx^2\right)  \ \textrm{and} \ 
\left(I_2, h_2=-dy^2, \bar h_2 = -\frac{1}{Y(x)^2}dy^2\right). 
 $$
 We assume that $X(x)>Y(y)$ for all $(x,y)$. 
 The  corresponding tensors $L_1$ and  $L_2$ (we view them as $1\times 1$-matrices)  and their characteristic polynomials are 
 $$
 L_1= (X(x))\ ; \ \ L_2=(Y(y)) \ ; \ \ \chi_1(t)= t - X(x) \ ; \  \   \chi_2(t)= t - Y(y).
 $$ 
 We see that the metrics $h_1$, $h_2$ satisfy the assumptions in Theorem \ref{thm3}.
 Plugging these data in the formulas \eqref{matg}, we obtain geodesically equivalent  metrics  $g$ and $\bar g$ 
 given by the matrices 
 $$
 g =\begin{pmatrix} X(x)- Y(y) & \\ & X(x)- Y(y)\end{pmatrix} \ ,  \ \  \bar g =
 \begin{pmatrix} \tfrac{X(x)- Y(y)}{X(x)^2Y(y)} & \\ & \tfrac{X(x)- Y(y)}{X(x)Y(y)^2}.\end{pmatrix}
 $$
  We see that these metrics are precisely the   Dini metrics  \eqref{dini}.
 For further use let us note that the tensor \eqref{L} for these metrics is given by $ L= \begin{pmatrix} X(x) & \\ &  Y(y)\end{pmatrix}.$
  \end{Ex}

  \begin{Ex}[Levi-Civita metrics (\ref{LC1},\ref{LC2})   follow from  splitting-gluing constructions.]
We take 4 pairs of geodesically equivalent metrics on the interval $I$.
\begin{equation} \label{met}
g_1= dx_1^2 \sim \bar g_1 = \tfrac{1}{X_1(x_1)^2}dx_1^2\ ; \ \ 
g_2= -dx_2^2 \sim \bar g_2 = -\tfrac{1}{X_2(x_2)^2}dx_2^2\ ; $$ $$
 g_3= dx_3^2 \sim \bar g_3 = \tfrac{1}{X_3(x_3)^2}dx_3^2\ ; \ \ 
g_4= -dx_4^2 \sim \bar g_4 = -\tfrac{1}{X_4(x_4)^2}dx_4^4. \end{equation}  We assume that for $i\ne j$ 
$X_i(x_i)\ne X_j(x_j)$ for all $x_i,x_j\in I$. 

Gluing  $(I, g_1, \bar g_1)$  and $(I, g_2, \bar g_2)$, 
($(I, g_3, \bar g_3)$  and $(I, g_4, \bar g_4)$,respectively)   we obtain two pairs of  geodesically equivalent metrics
(we denote them by $h_1\sim \bar h_1$ \ ($h_2 \sim \bar h_2$, respectively))  on the two-dimensional disk  
$U^2=I\times I$. These metrics and the corresponding tensors \eqref{L} were essentially constructed in Example \ref{diniformulas} and are given by   matrices 
$$
h_1 =\begin{pmatrix} X_1(x_1)- X_2(x_2) & \\ & X_1(x_1)- X_2(x_2)\end{pmatrix} \ \ \sim   \ \  \bar h_1 =
 \begin{pmatrix} \tfrac{X_1(x_1)- X_2(x_2)}{X_1(x_1)^2X_2(x_2)} & \\ & \tfrac{X_1(x_1)- X_2(x_2)}{X_1(x_1)X_2(x_2)^2}\end{pmatrix} \ ,  $$
  $$
h_2=    \begin{pmatrix} X_3(x_3)- X_4(x_4) & \\ & X_3(x_3)- X_4(x_4)\end{pmatrix} \ \sim   \  \bar h_2=
 \begin{pmatrix} \tfrac{X_3(x_3)- X_4(x_4)}{X_3(x_3)^2X_4(x_4)} & \\ & \tfrac{X_3(x_3)- X_4(x_4)}{X_3(x_4)X_4(x_4)^2}\end{pmatrix}\ ,  $$ $$L_1=L(h_1, \bar h_1)=  \begin{pmatrix} X_1(x_1) & \\ &  X_2(x_2)\end{pmatrix} \ , \ \
 L_2=L(h_2, \bar h_2)=  \begin{pmatrix} X_3(x_3) & \\ &  X_4(x_4)\end{pmatrix}.
$$
 We see that the metrics $h_1$, $h_2$ satisfy the assumptions in Theorem \ref{thm3}.
Gluing these metrics, we obtain the metrics  (\ref{LC1},\ref{LC2}). 

\end{Ex} 

\begin{Rem} By changing the sign of the metrics \eqref{met} we can make  geodesically equivalent 
metrics $g\sim \bar g$  of  arbitrary signature. \end{Rem}  

\begin{Ex}[General Levi-Civita metrics]
We take $m$ building blocks: the first $r$ building blocks are  1-dimensional,  
  and the last $m-r$ building blocks $h_{r+1}\sim \bar h_{r+1},...,h_{m}\sim \bar h_{m}$
    have dimensions   $k_i\ge 2$, $i=r+1, ..., m-r$.  For cosmetic
   reasons we think that  the first $r$ building blocks are 
 \begin{equation} \label{pm}
 (U_i^1, h_i=  \pm {dx_i}^2 , \bar h_i=  \pm  \frac{1}{X_i(x_i)^2} {dx_i}^2) \ , \ \ {i=1,...,r},\end{equation}
the sign $\pm$ in $h_i$ and $\bar h_i$  is the same for each $i$, but may be different for different $i$'s.   
The last $m-r$ building blocks are 
 $$
\left(U_i^{k_i}, h_i= \sum_{\alpha_i, \beta_i=1}^{k_i}(h_i(x_i))_{\alpha_i \beta_i}dx_i^{\alpha_i}dx_i^{\beta_i}, \bar h_i=  \frac{1}{X_i^{k+1}}\sum_{\alpha_i, \beta_i=1}^{k_i}(h_i(x_i))_{\alpha_i \beta_i} dx_i^{\alpha_i}dx_i^{\beta_i}\right) \ , \ \ {i=r+1,...,m}.$$

 Here the functions $X_i$ are constant for $i>r$ and depend only on the corresponding variable $x_i$ for $i\le r$. As above, we assume that $\textrm{Image}(X_i)\cap \textrm{Image}(X_j)=\varnothing$ for $i\ne j$. 
 The metrics $ h_i$, $i=r+1,...,m$ can be arbitrary, but  their  entries  $(h_i)_{\alpha_i\beta_i}$  must depend on the coordinates $x_i= (x_i^1,...,x_i^{k_i})$ only.

Inductively applying the gluing procedure,  we obtain for $g$ and $\bar g$ the following  form:
 
\begin{equation} \label{LCM}
\begin{array}{cccc} g &  = &  \sum_{i=1}^rP_i{dx_i}^2  & + \sum_{i=r+1}^m \left[P_i \sum_{\alpha_i, \beta_i=1}^{k_i}(h_i(x_i))_{\alpha_i \beta_i}dx_i^{\alpha_i}dx_i^{\beta_i}\right],  \\ 
{\bar g}  &  = &   \sum_{i=1}^r P_i \rho_i {dx_i}^2 & + \sum_{i=r+1}^m \left[ P_i \rho_i\sum_{\alpha_i, \beta_i=1}^{k_i}(h_i(x_i))_{\alpha_i \beta_i} dx_i^{\alpha_i}dx_i^{\beta_i}\right], 
\end{array}\end{equation}
 where 
 \begin{equation} \label{P_i}
 P_i:= \pm \prod_{j\ne i} (X_i- X_j), \  \  \ \rho_i:=  \frac{1}{X_i \, \prod_{\alpha} X_\alpha}. 
 \end{equation}
(the signs $\pm$  in \eqref{P_i} depend on the choice of  the signs $\pm$ in \eqref{pm} and can be arbitrary). 
This is precisely   Levi-Civita's normal form for geodesically equivalent (Riemannian)
metrics from \cite{Levi-Civita}. 

Now, since every pair of geodesically equivalent 
metrics (in a neighborhood of almost every point)
   can be obtained by a gluing construction, and since in the Riemannian signature only the blocks used above can be used, every Riemannian geodesically equivalent metrics have the form \eqref{LCM} in a certain coordinate system. This is the famous Levi-Civita's
    Theorem  from \cite{Levi-Civita}.
  \end{Ex}

    Note, than the Lorentz signature  of $g$ and $\bar g$ does not allow the tensor $L$ to have complex eigenvalues  
     of algebraic multiplicity greater than one. Similarly, it does not allow the tensor $L$ to have a Jordan block of dimension 4, or two Jordan blocks. Thus,   in order to obtain the description of nonpropotional
      4-dimensional geodesically      
equivalent metrics of Lorentz signature, one needs    the building blocks of dimensions $1,2,3$ only. In dimension 1, only one buiding block, namely the one from Example \ref{bb1}, is possible.

Geodesically equivalent metrics such that the tensor $L$ has  the 2-dimensional  
 Jordan-block structure  $$\begin{pmatrix} 
\lambda &1 \\ 
& \lambda \end{pmatrix}\ , \ \ 
 \begin{pmatrix} 
\lambda & \\ 
& \lambda \end{pmatrix}\ , \ \
\begin{pmatrix} 
\alpha &\beta& \\ 
-\beta & \alpha  
\end{pmatrix}.$$
 were described in  Example \ref{bb2}. 
  For the Jordan-block structure \begin{equation}\label{petrovcase}\begin{pmatrix} 
\lambda &1& \\ 
& \lambda &1 \\ &&
\lambda\end{pmatrix},\end{equation} 
 the description of the 
 metrics  follows from   Petrov \cite{Petrov49}:  the metrics are given by  
\begin{equation} \label{bolsinov}  
\begin{array}{ccl} g & =& \left(4\, x_{{2}}\left( {\frac {d}{dx_{{3}}}}\lambda
 \left( x_{{3}} \right)  \right) +2\,\right){ dx}_{{1
}}{ dx}_{{3}} +{{ dx}_{{2}}}^{2}\\ &+&2\, x_{{1}} \left( {\frac 
{d}{dx_{{3}}}}\lambda \left( x_{{3}} \right)  \right)  { dx}_{{2}}{ dx}_{{3}}+ {x_{{1}}}^{2}
\left( {\frac {d}{dx_{{3}}}}\lambda \left( x_{{3}}
 \right)  \right) ^{2}{{ dx
}_{{3}}}^{2} ,\\
\bar g & =& 
\frac{1}{  \lambda \left( x_{{3}} \right)^{6}} 
\Big[\left.\left(4\, x_{{2}}  \lambda \left( 
x_{{3}} \right)^{2}\left( {\frac {d}{dx_{{3}}}}
\lambda \left( x_{{3}} \right)  \right) +2\,  
\lambda \left( x_{{3}} \right)^{2}\right){ dx}_{{1}}{ dx}_{{3}}
+  \lambda \left( x_{{3}} \right)^{2}{{ dx}_{{2}}}^{2}
 \right. \\ 
  &-&\left. \left(4\, x_{{2}}\lambda \left( x_{{3}} \right)\left( {\frac {d}{dx_{{3}}}}\lambda \left( x_{{3}}
 \right)  \right)  +2\lambda \left( x_{{3}} \right) -2\,x_{{1}}  \lambda \left( x_{{3}} \right)^2\left( {\frac {d}{dx_{{3}}}}\lambda \left( x_{{3}}
 \right) \right)  \right) {dx}_{{2}}
{ dx}_{{3}} 
 \right. \\
  &+& \left. \left( 4\, {x_{{2}}}^{2} \left( {\frac {d}{dx_{{3}}}}
\lambda \left( x_{{3}} \right)  \right)^{2}+4\,x_{{2}} \left( {\frac {d}{dx_{{3}}}}\lambda \left( x_{{3}}
 \right)  \right) - 4\, x_{{1}}x_{{2}}
\lambda \left( x_{{3}} \right)\left( {\frac {d}{dx
_{{3}}}}\lambda \left( x_{{3}} \right)  \right)^{2}\right){{ dx}_{{3}}}^{2}  \right. 
\\ &+& \left.  \left(1+ 
 {x_{{1}}}^{2}  \lambda \left( x_{{3}} \right)^{2}\left( {\frac {d}{dx_{{3}}}}\lambda \left( x_{{3}} \right)  \right)^{2} -2\,x_{{1}}\lambda \left( x_{{3}} \right)\left( {\frac {d}{dx_{{3}}}}\lambda \left( x_{{3}} \right) 
 \right)\right)  {{ dx}_{{3}}}^{2}\Big]\right.
\end{array}
 \end{equation}

The corresponding $L$ is given by the matrix $$\left[ \begin {array}{ccc} \lambda \left( x_{{3}} \right) &1& \left( 
{\frac {d}{dx_{{3}}}}\lambda \left( x_{{3}} \right)  \right) x_{{1}}
\\ \noalign{\medskip}0&\lambda \left( x_{{3}} \right) &2\, \left( {
\frac {d}{dx_{{3}}}}\lambda \left( x_{{3}} \right)  \right) x_{{2}}+1
\\ \noalign{\medskip}0&0&\lambda \left( x_{{3}} \right) \end {array}
 \right]. $$

\begin{Rem}    
  Actually, the formulas  \eqref{bolsinov}   are slightly more  general than that of \cite{Petrov49}. They are equivalent to the formulas from \cite{Petrov49} (modulo a coordinate transformation) 
   at the  points such that $d\lambda\ne 0$. The formulas \cite{Petrov49} were obtained together  with A. Bolsinov; they can be generalized for  every dimension. We will publish this result elsewhere.
  \end{Rem} 
  
  As it follows from \cite[Lemma 6]{einstein}, if $L$ has the Jordan-form  $\begin{pmatrix} 
\lambda &1& \\ 
& \lambda & \\ &&
\lambda\end{pmatrix}$, the eigenvalue $\lambda$ is constant, and the metrics are affinely 
 equivalent (i.e., Levi-Civita connections of $g$ and $\bar g$ coincide). Affinely equivalent metrics whose tensor $L$ has this form were essentially described by Eisenhart in \cite{eisenhart38},   see also \cite[Theorem 1]{solodovnikov1959}. From their description it follows, that, 
  in a certain coordinate system,  geodesically equivalent metrics $g\sim \bar g$ are given   by  
\begin{equation} \label{metricseisenhart} \begin{array}{cl} g&=  2\,{ dx}_{{3}}{ dx}_{{1}}+
 { h}\left( x_{{2}},x_{{3}} \right)_{11}{{ dx}_{{2}}}^{2}+2\, { h(x_2,x_3)_{12}} { dx}_{{2}}{ dx}_{{3}}
  +{ h(x_2,x_3)_{22}}{{ dx}_{{3}}}^{2},  \\  
  \bar g &= 2\,\alpha\,{ dx}_{{3}}{ dx}_{{1}}+\alpha\,{
 h\left( x_{{2}},x_{{3}} \right)_{11}}{{ dx}_{{2}}}^{2}  
 +2\,\alpha\,{ h\left( x_{{2}},x_{{3}} \right)_{12}} { dx}_{{2}}{ dx}_{{3
}} +\beta {{ dx}_{{3}}}^{
2}+\alpha\,{ h\left( x_{{2}},x_{{3}}
 \right)_{22}}{{ dx}_{{3}}}^{2}, \end{array}
 \end{equation} where $\alpha$ and $\beta$ are constants.

  Now, the metrics $g,\bar g$ such that   $L= \begin{pmatrix} 
\lambda &  &\\ 
& \lambda & \\ & & \lambda\end{pmatrix} $ are conformally equivalent. By 
 by the classical result of Weyl \cite{Weyl}, they are proportional (i.e., $ \bar g =\const\cdot  g$).

  Thus, we have described all building blocks  that can be used in constructing  
  metrics of  Lorentz signature; Theorem \ref{thm3} gives us the construction.  Let us count the number of cases in dimension 4:  we can represent $4$ as the sum of natural numbers by 4 different ways: 
  
  \begin{tabular}{c|c|c} \hline
  Dim of blocks & Description of blocks & \# of cases\\ \hline  
  1+1+1+1 & \begin{minipage}{.6\textwidth}All  building  blocks are as in Example \ref{bb1}, and  $g\sim \bar g$  are  essentially (\ref{LC1},\ref{LC2}) with the changed sign of $dx_1^2$\end{minipage}&  1\\ \hline  
 1+1+2  &  \begin{minipage}{.6\textwidth}The first two building blocks  are  as in Example \ref{bb1}, the third is as in Example \ref{bb2}\end{minipage}&3\\ \hline   
 2+2    & \begin{minipage}{.6\textwidth}Both building  blocks are as in Example \ref{bb2}; at least one of them is trivial \end{minipage}& 3 \\
 \hline  
 1+3 & \begin{minipage}{.6\textwidth} The first building  block as is Example \ref{bb1}, the second   is as in  \eqref{bolsinov}, as in \eqref{metricseisenhart}, or trivial \end{minipage}& 3 \\ \hline 
  \end{tabular}

  \begin{Rem} The general schema also works in higher dimensions, but in this case there is the following  essential difficulty (and this is the only difficulty):  up to our knowledge, for dimensions $n-1\ge 5$,  
  there is no description of all pairs of $(g, L)$ such that $g$ has Lorentz signature and  $L$ is an (1,1)-selfadjoint tensor such that it is covariantly constant, and such that the Jordan normal form of $L$ is 
  \begin{equation}\label{llll}
  \begin{pmatrix} 
\lambda &1& &&&\\ 
& \lambda &1 &&& \\
&& \lambda&0&&\\
&&&\ddots&\ddots&\\
 &&&&\lambda & 0\\
&&&&&\lambda\end{pmatrix} \end{equation}   
  In  dimension $n  = 4$, since $n-1=3$,    the Jordan normal form \eqref{llll} coincides with \eqref{petrovcase}, and the 
    local description follows from  \cite{Petrov49}.  In dimension $n=5$ we have $n-1=4$ and one can obtain the local description (we will not do it in the present paper) combining the results of \cite{eisenhart38,solodovnikov1959} with the algebraic description of possible holonomy groups of   4-dimensional metrics of Lorentz signature (see e.g. \cite{ new2, new3}). 
\end{Rem} 
  
  \paragraph{\bf Acknowledgement.}
This work benefited from discussions with A. Bolsinov, G. Gibbons, D. Giulini,   V. Kiosak, P. Nurowski,  and A. Wipf. I thank the anonimous referee and G. Hall for valuable suggestions and finding misprints.   During the work on this paper, the author was partially supported by  
Deutsche Forschungsgemeinschaft (SPP 1154 and GK 1523) and FSU Jena.


\begin{thebibliography}{99}

\bibitem{Aminova} A. V.  Aminova, {\it Pseudo-Riemannian manifolds with general geodesics,}   Russian Math. Surveys  {\bf 48}(1993),  no. 2, 105--160.

\weg{\bibitem{Aminova2}
A. V. Aminova,
\emph{Projective transformations of pseudo-Riemannian manifolds. Geometry, 9.}
 J. Math. Sci. (N. Y.) \textbf{113}(2003),  no. 3, 367--470.}

\bibitem{Beltrami}
E. Beltrami,
\emph{Risoluzione del problema: riportare i punti di una
superficie sopra un piano in modo che le linee geodetiche vengano rappresentate da linee rette,
}
Ann. Mat., \textbf{1}(1865),  no. 7, 185--204.

\bibitem{Beltrami2} E. Beltrami, \emph{Saggio di interpetrazione della geometria non-euclidea,}  
Giornale di matematiche, vol. {\bf VI}(1868). 


\bibitem{Beltrami3} E. Beltrami, \emph{Teoria fondamentale degli spazii di curvatura costante,}
 Annali. di Mat., ser II {\bf 2}(1968),  232--255.

\bibitem{BMF} A. V. Bolsinov, V. S.  Matveev,  A. T.  Fomenko,
\emph{  Two-dimensional Riemannian metrics with an integrable
geodesic flow. Local and global geometries,} Sb. Math. {\bf
189}(1998),  no. 9-10, 1441--1466.

\bibitem{pucacco} A. V. Bolsinov, V. S.  Matveev, G. Pucacco, {\it 
Normal forms for pseudo-Riemannian 2-dimensional metrics whose geodesic flows admit integrals quadratic in momenta,} J. Geom. Phys. {\bf 59}(2009), no. 7,   1048--1062.  	arXiv:math.DG/0803.0289v2 


\bibitem{splitting} A. V. Bolsinov, V. S.  Matveev,
{\it Splitting and gluing lemmas for geodesically equivalent pseudo-Riemannian metrics,} accepted to  Transactions of the American Mathematical Society. 	arXiv:math.DG/0904.0535.

\bibitem{Brozos}  M. Brozos-V\'azquez, P. Gilkey,   H.  Kang, S. Nikcevic, G.  Weingart, {\it 
Geometric realizations of curvature models by manifolds with constant scalar curvature,} 
Differential Geom. Appl. {\bf 27}(2009), no. 6, 696--701.


\bibitem{bryant} R. L. Bryant, G. Manno, V. S. Matveev, \emph{A
    solution of a problem of Sophus Lie: Normal forms of 2-dim metrics
    admitting two projective vector fields}, Math. Ann. {\bf 340}, no. 2, 437--463, 2008. 	arXiv:0705.3592
    
\bibitem{BDE} R. L. Bryant, M. Dunajski, M. Eastwood,
  \emph{Metrisability of two-dimensional projective structures},
  J. Diff. Geom. {\bf 83}(2009), no. 3, 465--499.

\bibitem{Darboux}
G. Darboux,
\emph{Le\c{c}ons sur la th\'eorie g\'en\'erale des surfaces},
Vol. III, Chelsea Publishing, 1896.

\bibitem{eastwood} 
M. Eastwood,  V. S. Matveev,  \emph{ Metric connections in projective differential geometry,}
 Symmetries and Overdetermined Systems of Partial Differential Equations (Minneapolis, MN, 2006), 339--351,
 IMA Vol. Math. Appl.,    {\bf
144}(2007),   Springer, New York.  	arXiv:0806.3998.

\bibitem{Ehlers} J. Ehlers, F. A.  E. Pirani, A. Schild,  \emph{The geometry of free fall and light propagation.}  General relativity (papers in honour of J. L. Synge), Clarendon Press, Oxford, 1972,  
 63--84. 
 \bibitem{eisenhart23} L. P. Eisenhart, \emph{
The geometry of paths and general relativity,}
Ann. of Math. (2) {\bf 24}(1923), no. 4, 367--392. 

\bibitem{eisenhart38} L. P. Eisenhart, \emph{
Fields of parallel vectors in Riemannian space.}
Ann. of Math. (2) {\bf 39}(1938), no. 2, 316--321, see also  

\bibitem{Gibbons} G. W. Gibbons, C. M. Warnick, \emph{Dark Energy and Projective Symmetry,} Physics Letters B {\bf 688}(2010),  337--340.	arXiv:1003.3845. 


\bibitem{Gilkey2001} P. Gilkey, \emph{ Geometric properties of natural operators defined by the Riemann curvature tensor,}  World Scientific Publishing Co.,  2001, viii+306 pp. 

\bibitem{Golikov} V. I. Golikov, {\it Geodesic mappings of gravitational fields of general type,} 
 Trudy Sem. Vektor. Tenzor. Anal., {\bf 12}(1963)  79--129.
 
 \bibitem{Hall83} G. S. Hall, {\it Curvature collineations and the determination of the metric from the curvature in general relativity,}  Gen. Relativity Gravitation {\bf 15}(1983), no. 6, 581--589. 
 
 {\bibitem{rendall} G. S. Hall, A. D. Rendall,  {\em Uniqueness of the metric from the Weyl and energy-momentum tensors,}  J. Math. Phys. {\bf  28}(1987),  no. 8, 1837--1839.}
 
 \bibitem{book} G. S. Hall, {\em  Symmetries and curvature structure in general relativity, } World Scientific Lecture Notes in Physics, 46. World Scientific Publishing Co., Inc., River Edge, NJ, 2004. x+430 pp.
 
 \bibitem{Hall2007} G. S. Hall, D. P.  Lonie,
 \emph{
The principle of equivalence and  projective structure in spacetimes,}
Classical Quantum Gravity {\bf 24}(2007),  14, 3617--3636. 
 
 \bibitem{Hall2008} G. S. Hall,  D. P. Lonie, \emph{  The principle of equivalence and cosmological metrics}, J. Math.
Phys. {\bf 49}(2008), 022502. 
 
 \bibitem{Hall2010} G. S. Hall, D. P. Lonie,  \emph{
Projective equivalence of Einstein spaces in general relativity,} 
Classical Quantum Gravity  {\bf 26}(2009), no. 12, 125009, 10 pp.
 
 
 
\bibitem{new2} G. S. Hall, D. P. Lonie,  \emph{ Holonomy and projective equivalence in 4-dimensional Lorentz manifolds,} SIGMA {\bf 5}(2009), Paper 066, 23 pp.

\bibitem{new3} G. S. Hall, D. P. Lonie,  \emph{Projective structure and holonomy in four-dimensional Lorentz manifolds,}  Journal of Geometry and Physics
{\bf 61}(2011) no.  2,  381--399 

\bibitem{einstein} 
V. Kiosak,  V. S. Matveev,  \emph{ Complete Einstein metrics are geodesically rigid,}    Comm. Math. Phys. {\bf 289}(2009), no. 1, 383-400.  
  	arXiv:0806.3169.
  	
\bibitem{KioMat2010} V. Kiosak, V. S. Matveev, \emph{Proof Of The
    Projective Lichnerowicz Conjecture For Pseudo-Riemannian Metrics
    With Degree Of Mobility Greater Than Two}, Comm. Mat. Phys. {\bf  297},
  no. 2, 401--426, 2010  	
  	
 \bibitem{solodovnikov1959} G. I. Kruckovic,  A. S.  Solodovnikov, \emph{Constant symmetric tensors in Riemannian spaces,}  Izv. Vys. Uchebn. Zaved. Matematika {\bf 1959} no. 3(10), 147--158.
 
\bibitem{Kruchkovich}  G. I. Kruchkovich, {\it Equations of semireducibility and geodesic correspondence of Lorentz spaces,}  Trudy Vsecsoyuz. Zaochn. Energet. Inst., {\bf 24}(1963), 74--87.



\bibitem{kruglikov2008} B. Kruglikov, \emph{ Invariant characterization of Liouville metrics and polynomial integrals,} J.
Geom. Phys. {\bf 58}(2008), no. 8, 979--995. arXiv:0709.0423


\bibitem{lagrange} J.-L.  Lagrange, \emph{ Sur la construction des cartes g\'eographiques,} Nov\'eaux M\'emoires de l'Acad\'emie des Sciences et Bell-Lettres de Berlin, 1779. 

 \bibitem{Levi-Civita}
 T. Levi-Civita, {\it Sulle trasformazioni delle equazioni
 dinamiche}, Ann. di Mat., serie $2^a$, {\bf 24}(1896), 255--300.
 
 \bibitem{Lie} 
S. Lie, 
\emph{Untersuchungen \"uber geod\"atische Kurven},
Math. Ann. \textbf{20} (1882); Sophus Lie Gesammelte Abhandlungen, 
Band 2, erster Teil, 267--374. Teubner, Leipzig, 1935.

 \bibitem{Liouville}
R. Liouville, 
\emph{Sur les invariants de certaines \'equations diff\'erentielles
et sur leurs applications},
Journal de l'\'Ecole Polytechnique \textbf{59} (1889), 7--76.

  \bibitem{bifurcations}
V. S. Matveev,  {\it On projectively equivalent metrics near points of bifurcation,}  In the book ``Topological methods in the theory of integrable systems''(Eds.: Bolsinov A.V., Fomenko A.T., Oshemkov A.A.), Camb. Sci. Publ., 2006,  pp. 214 -- 240, arXiv:0809.3602.  
 
 \bibitem{archive}
V. S. Matveev,
\emph{Proof of projective Lichnerowicz-Obata conjecture},
 J. Diff. Geom.  {\bf 75}(2007),  459--502,  arXiv:math/0407337
 
 \bibitem{mcintosh}  C. B. G. McIntosh,  W. D. Halford, \emph{
Determination of the metric tensor from components of the Riemann tensor,}
J. Phys. A {\bf 14}(1981), no. 9, 2331--2338. 
 
 \bibitem{Nurowski2}	P. Nurowski, \emph{Projective vs metric structures,}  J. Geom.  Phys., accepted,   arXiv:1003.1469 
  
  \bibitem{Nurowski}	P. Nurowski, \emph{Is dark energy meaningless?}
  Rendiconti del Seminario Matematico
Universita  e Politecnico di Torino, {\bf 68} (2010), no. 4,  361--367,  arXiv:1003.1503.   




\bibitem{Petrov49} A. Z. Petrov, \emph{Geodesic mappings of Riemannian spaces of an indefinite metric (Russian)}, Uchen. Zap. Kazan. Univ., {\bf 109}(1949),   no. 3, 7--36.

\bibitem{Petrov1} A. Z. Petrov, \emph{  On a geodesic representation of Einstein spaces (Russian),}   Izv. Vys. Ucebn. Zaved. Matematika    {\bf 21}(1961) no. 2, 130--136.

\bibitem{sinjukov54} N. S.  Sinjukov,
 {\it
On geodesic mappings of Riemannian spaces onto symmetric Riemannian spaces, } 
Dokl. Akad. Nauk SSSR (N.S.) {\bf 98}(1954), 21--23. 

\bibitem{sinjukov} N. S.  Sinjukov,
 {\it Geodesic mappings of Riemannian spaces},  (in Russian)
``Nauka'', Moscow, 1979. 

\bibitem{thomas} T. Thomas, \emph{
On the projective theory of two dimensional Riemann spaces,}
Proc. Nat. Acad. Sci. U. S.  A. {\bf  31}(1945) 259--261.

\bibitem{veblen23} O. Veblen, T.  Thomas, \emph{
The geometry of paths,}
Trans. Amer. Math. Soc. {\bf 2}(1923), no. 4, 551--608.

\bibitem{veblen26} O. Veblen, J.  Thomas, \emph{ Projective invariants of affine geometry of paths,} Ann. of Math. (2) {\bf 27}(1926), no. 3, 279--296.

\bibitem{Weyl}  H. Weyl, {\it Zur Infinitisimalgeometrie: Einordnung der projektiven
und der  konformen Auffasung,} Nachrichten von der K. Gesellschaft
der Wissenschaften zu G\"ottingen, Mathematisch-Physikalische
Klasse, 1921;
 ``Selecta Hermann Weyl'', Birkh\"auser Verlag,
   Basel und Stuttgart,
1956. 

\end{thebibliography}
\end{document}